\documentclass[11pt,a4paper]{article}
\usepackage[T1]{fontenc}
\usepackage[utf8]{inputenc}
\usepackage{lmodern}
\usepackage[margin=27mm]{geometry}
\usepackage{amsmath,amssymb,amsthm,mathtools}
\usepackage{microtype,enumitem,needspace}
\usepackage{booktabs,tabularx}
\usepackage{tikz}
\usepackage[numbers,sort&compress]{natbib}
\usepackage[hidelinks]{hyperref}
\usepackage{doi}
\usepackage[nameinlink,capitalise,noabbrev]{cleveref}
\usepackage{etoolbox}
\newtheorem{theorem}{Theorem}[section]
\newtheorem{lemma}[theorem]{Lemma}
\newtheorem{corollary}[theorem]{Corollary}
\newtheorem{proposition}[theorem]{Proposition}
\theoremstyle{definition}

\newtheorem{problem}{Problem}
\crefname{problem}{Problem}{Problems}
\Crefname{problem}{Problem}{Problems}
\theoremstyle{remark}

\AtBeginEnvironment{lemma}{\crefalias{theorem}{lemma}}
\AtBeginEnvironment{corollary}{\crefalias{theorem}{corollary}}
\AtBeginEnvironment{proposition}{\crefalias{theorem}{proposition}}
\AtBeginEnvironment{remark}{\crefalias{theorem}{remark}}
\AtBeginEnvironment{definition}{\crefalias{theorem}{definition}}
\DeclareMathOperator{\comp}{comp}
\newcommand{\E}{\mathbb E}
\newcommand{\Prob}{\mathbb P}
\newcommand{\ind}{\mathbf 1}
\newcommand{\dd}{\,\mathrm d}
\newcommand{\join}{\mathbin\vee}
\newcommand{\norm}[1]{\lVert#1\rVert}
\newcommand{\abs}[1]{\lvert#1\rvert}
\newcommand{\inner}[2]{\langle#1,#2\rangle}
\numberwithin{equation}{section}
\allowdisplaybreaks[1]
\title{Goodman-type lower bounds for graph homomorphism densities}
\author{Taeyoung Kim\\[0.4em]
  \small School of Computing, KAIST, Republic of Korea\\
  \small School of Computational Sciences, KIAS, Republic of Korea}
\date{}
\begin{document}
\maketitle

\begin{abstract}
For a finite graph $H$ with $v(H)$ vertices, chromatic number
$\chi(H)$, and chromatic polynomial $P_H$, define
$\Phi_H(p)=(1-p)^{v(H)}P_H((1-p)^{-1})$. This polynomial interpolates
the homomorphism densities of $H$ in balanced complete multipartite
graphons. We study the lower bound $t(H,W)\ge\Phi_H(p)$ on the
homomorphism density of $H$ in graphons
$W$ of edge density $p>1-1/(\chi(H)-1)$, which generalises Goodman's
triangle bound. We prove the inequality for chordal graphs whose
maximal cliques have equal size, for trivially perfect graphs, and
for joins of a clique and an odd cycle. The proofs use entropy along
clique trees and a convexity criterion for adding a universal vertex.
We also prove that a nonbipartite graph satisfies the inequality at
all sufficiently high edge densities if and only if its girth is odd.
For a graph on $n$ vertices whose girth is odd, the bound holds whenever
$p>1-(72n^{12})^{-1}$. The proof combines Whitney's forest expansion
with a spectral comparison of cycle densities. Finally, among the
$117$ nonbipartite graphs on at most six vertices with no isolated
vertices or isolated edges, exactly $94$ satisfy the bound. The main
results are formally verified in Lean~4.
\end{abstract}

\noindent\textit{Keywords:} graphon, homomorphism density, chromatic
polynomial, Goodman's inequality, chordal graph, Tur\'an graphon.

\smallskip
\noindent\textit{Mathematics Subject Classification (2020):}
05C35 (primary); 05C15, 05C31 (secondary).

\section{Introduction}\label{sec:introduction}

Goodman's inequality~\cite{Goodman1959} gives the lower bound
$p(2p-1)$ for the triangle homomorphism density of a graphon of edge
density $p$. Equality holds for balanced complete $k$-partite
graphons, whose edge densities are $1-1/k$. The exact minimum clique
density at a given edge density was determined by
Razborov~\cite{Razborov2008} for triangles, by
Nikiforov~\cite{Nikiforov2011} for $K_4$, and by
Reiher~\cite{Reiher2016} for all cliques. In this paper we ask which
graphs $H$ satisfy the polynomial lower bound obtained by interpolating
their densities in balanced complete multipartite graphons.

All graphs are finite and simple. For a graph $H$, write $V(H)$ and
$E(H)$ for its vertex and edge sets, $v(H)$ and $e(H)$ for their
cardinalities, $d_H(v)$ for the degree of a vertex $v$, and $\chi(H)$
for the chromatic number. Let $K_s$, $P_s$, and $C_s$ denote the
complete graph, path, and cycle, respectively, each on $s$ vertices; $K_0$
is the empty graph. We use $\sqcup$ for disjoint union.
A \emph{graphon} is a symmetric measurable function
$W\colon[0,1]^2\to[0,1]$; see Lov\'asz~\cite{Lovasz2012}.
The homomorphism densities of $W$ are defined by
\[
 t(H,W)=\int_{[0,1]^{V(H)}}\prod_{uv\in E(H)}W(x_u,x_v)
                  \prod_{v\in V(H)}\dd x_v,
 \qquad p=t(K_2,W).
\]
Let $P_H$ be the chromatic polynomial of $H$ and define
\begin{equation}\label{eq:target}
 \Phi_H(p)=(1-p)^{v(H)}P_H\!\left(\frac1{1-p}\right).
\end{equation}
The expression in \eqref{eq:target}, initially defined for $p<1$,
extends to a polynomial with $\Phi_H(1)=1$, since $P_H$ is monic of
degree $v(H)$. For $\chi(H)\ge2$, we say that $H$ satisfies the
\emph{Goodman-type lower bound} if every graphon $W$ satisfies
\begin{equation}\label{eq:goodman}
 t(H,W)\ge\Phi_H(p)\qquad\text{whenever}\qquad
 p>\underline{p}_H:=1-\frac1{\chi(H)-1}.
\end{equation}
For $H=K_3$, we have $\Phi_H(p)=p(2p-1)$.

For an integer $k\ge1$, the \emph{Tur\'an graphon} $T_k$ is obtained
by partitioning $[0,1]$ into $k$ sets of equal size and assigning
value zero within each part and value one between distinct parts.
It satisfies
\begin{equation}\label{eq:turan-values}
 t(K_2,T_k)=1-\frac1k,\qquad
 t(H,T_k)=\frac{P_H(k)}{k^{v(H)}}=\Phi_H(1-1/k).
\end{equation}
Indeed, $t(H,T_k)$ is the fraction of maps $V(H)\to\{1,\ldots,k\}$
that are proper colourings. We call the numbers $1-1/k$
\emph{Tur\'an densities}, and we say that $T_k$ is \emph{optimal} for
$H$ if it minimises $t(H,W)$ among graphons $W$ of edge density
$1-1/k$, that is, if
\begin{equation}\label{eq:turan-optimal}
 t(K_2,W)=1-\frac1k\quad\Longrightarrow\quad
 t(H,W)\ge t(H,T_k).
\end{equation}
Thus \eqref{eq:goodman} implies that $T_k$ is optimal for $H$ for
every integer $k\ge\chi(H)$.

Im, Kim, and Yang~\cite[Section~8]{IKY} introduced the polynomial
\eqref{eq:target} and the problem of classifying the graphs that
satisfy \eqref{eq:goodman}. They also announced the chordal theorem
for maximal cliques of equal size (\cref{main:pure}), the results on clique
joins in \cref{sec:lifting}, and the obstruction for graphs of even
girth (\cref{evt:prop:even-obstruction}). We give complete proofs of
these results and establish the further criteria and classifications
stated below.

\subsection{Chordal graphs}

A graph is \emph{chordal} if it has no induced cycle of length at
least four. A \emph{clique tree} of a graph is a tree whose vertices
are its maximal cliques and in which the cliques containing any fixed
vertex form a connected subtree. Such a tree exists for every
connected chordal graph; see Gavril~\cite{Gavril} and the account of
clique trees by Blair and Peyton~\cite{BlairPeyton1993}.

\begin{theorem}[Chordal graphs with maximal cliques of equal size]\label{main:pure}
Let $H$ be a connected chordal graph whose maximal cliques all have
the same size $k\ge3$. Then $H$ satisfies the Goodman-type lower bound.
\end{theorem}

We deduce \cref{main:pure} from the following criterion, which also
allows maximal cliques of different sizes.

\begin{lemma}\label{main:profile}
Let $H$ be a connected chordal graph of clique number $k\ge3$.
Let $C_1,\ldots,C_m$ be its maximal cliques, and let $T$ be a clique
tree. Let $r_i=|C_i|$ and $s_e=|C_i\cap C_j|$ for $e=ij\in E(T)$,
and define
\begin{equation}\label{eq:bj}
 b_j=\#\{i:r_i\ge j\}-\#\{e\in E(T):s_e\ge j\},
 \qquad 1\le j\le k.
\end{equation}
If $b_1\le b_2\le\cdots\le b_k$, then $H$ satisfies the
Goodman-type lower bound.
\end{lemma}

As we show in \cref{sec:chordal}, $b_j$ is the multiplicity of $j-1$
as a root of $P_H$. The criterion therefore depends only on the
chromatic polynomial. If every maximal clique has size $k$, then
$b_j=m-\#\{e:s_e\ge j\}$ is nondecreasing, giving \cref{main:pure}.
The proof combines Lee's entropy construction along a
tree~\cite[Theorem~2.6]{Lee2021} with the Moon--Moser recurrence for
clique densities~\cite{MoonMoser1962}.

\subsection{Trivially perfect graphs and clique joins}

A graph is \emph{trivially perfect} if it has no induced $P_4$ or $C_4$;
see Golumbic~\cite{Golumbic1978}.
Equivalently, these graphs are generated from $K_1$ by disjoint union
and addition of a universal vertex, a vertex adjacent to every
existing vertex~\cite{DrangePilipczuk}.

\begin{theorem}[Trivially perfect graphs]\label{main:tp}
Every trivially perfect graph with at least one edge satisfies the
Goodman-type lower bound.
\end{theorem}

Write $F\join G$ for the join of $F$ and $G$, obtained from their
disjoint union by adding every edge between the two vertex sets.
We prove \cref{main:tp} by showing that addition of a universal
vertex preserves the bound under a convexity hypothesis. Specifically,
\cref{prop:lifting} shows that if $F$ satisfies
the Goodman-type lower bound and $\Phi_F$ is nonnegative,
nondecreasing, and convex on $[\underline{p}_F,1]$ with $\Phi_F(\underline{p}_F)=0$, then
$K_1\join F$ satisfies the bound, and $\Phi_{K_1\join F}$ has the same
properties on the corresponding interval. This also applies to odd
cycles, using the following result.

\begin{theorem}[Im, Kim, and Yang~\cite{IKY}]\label{main:odd}
For every odd integer $\ell\ge3$ and every graphon $W$ of edge density
$p>1/2$,
\begin{equation}\label{eq:odd-bound}
 t(C_\ell,W)\ge p^\ell-p(1-p)^{\ell-1}=\Phi_{C_\ell}(p).
\end{equation}
\end{theorem}

The range $p\ge2/3$ is proved in~\cite[Theorem~1.1]{IKY}, and the
extension to $p>1/2$ is announced in~\cite[Section~8]{IKY}. The
accompanying Lean development includes the proof on the full range;
see the paragraph on \hyperref[par:formalization]{formal verification}.
The polynomial $\Phi_{C_\ell}$ satisfies
the required convexity conditions on $[1/2,1]$, so repeated addition
of universal vertices gives the following consequence.

\begin{corollary}[Clique joins of odd cycles]\label{main:odd-join}
For every integer $s\ge0$ and every odd integer $\ell\ge3$, the graph
$K_s\join C_\ell$ satisfies the Goodman-type lower bound.
\end{corollary}

\subsection{High densities and girth}

We say that $H$ \emph{eventually satisfies} the Goodman-type lower
bound if there exists $p_H<1$
such that $t(H,W)\ge\Phi_H(p)$ holds for every graphon of edge
density $p>p_H$. For a graph containing a cycle, the \emph{girth} is
the length of a shortest cycle.

\Needspace{7\baselineskip}
\begin{theorem}[Eventual validity and girth]\label{main:eventual}
Let $H$ be a nonbipartite graph. The following are equivalent:
\begin{enumerate}[label=\textup{(\roman*)}]
 \item For all sufficiently large integers $k$, $T_k$ minimises
 $t(H,W)$ among graphons of edge density $1-1/k$.
 \item The graph $H$ eventually satisfies the Goodman-type lower
 bound.
 \item The girth of $H$ is odd.
\end{enumerate}
Moreover, if the girth of $H$ is odd, then $t(H,W)\ge\Phi_H(p)$ for
every graphon $W$ of edge density $p>p_H$, where
\begin{equation}\label{eq:eventual-polynomial-threshold}
 p_H=1-\frac{1}{72v(H)^{12}};
\end{equation}
in particular, $T_k$ is optimal for every integer $k>72v(H)^{12}$.
\end{theorem}

The proof in \cref{sec:eventual} writes $W = 1 - U$ and expands $t(H,1-U)$, then
groups the terms over forests using Whitney's broken-circuit theorem.
We estimate this sum through a spectral comparison of cycle densities
and bounds on the number of cycles in $H$. For even girth, the
constant graphon violates the bound at every density sufficiently
close to one.

\subsection{Graphs on at most six vertices}

Our last main result is a computer-assisted classification of small
graphs. We identify graphs by their numbers in Read and Wilson's
\emph{An Atlas of Graphs}~\cite{ReadWilson1998}.

\begin{theorem}[Classification on at most six vertices]\label{main:classification}
There are $117$ unlabelled nonbipartite graphs on at most six
vertices with no isolated vertices and no isolated edges. Of these, the $23$ graphs in
\cref{fig:catalogue-negative} violate the Goodman-type lower bound,
and the other $94$ satisfy it.
\end{theorem}

Removing an isolated vertex or an isolated edge from $H$ divides both
sides of \eqref{eq:goodman} by $1$ or by $p$, respectively, and does
not change $\chi(H)$. Hence a nonbipartite graph on at most six
vertices satisfies the Goodman-type lower bound if and only if the
graph obtained by deleting all its isolated vertices and isolated
edges does.

The proof combines general density inequalities, estimates for
individual graphs, and exact rational sum-of-squares certificates.
Explicit graphons give the counterexamples. Each of the $117$ cases
has been formally verified in Lean. In \cref{sec:catalogue}, we also
give density ranges on which the violating graphs satisfy the bound
and show that, for exactly $100$ of the $117$ graphs, $T_k$ is optimal
for every integer $k\ge\chi(H)$.

\subsection{Organisation}

\Cref{sec:chordal} proves the criterion for chordal graphs.
\Cref{sec:lifting} establishes the bound for trivially perfect graphs
and clique joins of odd cycles.
\Cref{sec:eventual} proves the characterisation by girth and the
quantitative threshold. \Cref{sec:catalogue} treats the classification
on at most six vertices and optimality at Tur\'an densities.
\Cref{sec:open-problems} discusses open problems, formal verification,
and the use of AI tools.

\section{A criterion for chordal graphs}\label{sec:chordal}

We first express the chromatic polynomial through a clique tree.
An entropy inequality then reduces the desired bound to a comparison
of successive clique densities.

\subsection{Clique-tree identities}

We use the clique tree and notation of \cref{main:profile}. For the
standard structural properties, see \cite{BlairPeyton1993}.
Root the tree at $C_1$ and order its cliques so that each parent
precedes its children. Set $S_1=\varnothing$. For each $C_i$ with parent
$C_{\pi(i)}$, set $S_i=C_i\cap C_{\pi(i)}$.
For every $i$, let $N_i=C_i\setminus S_i$.
The defining property of a clique tree gives
\[
 C_i\cap\bigcup_{h<i}C_h=S_i.
\]
The sets $N_i$ partition $V(H)$, and $|N_i|=r_i-|S_i|$.
Among the $\binom{r_i}{2}$ edges of $C_i$, precisely the
$\binom{|S_i|}{2}$ edges within $S_i$ have appeared in an earlier
clique. Counting each vertex and edge in the first clique containing
it therefore gives
\begin{equation}\label{eq:tree-counts}
 \begin{aligned}
 v(H)&=\sum_{i=1}^m(r_i-|S_i|)
       =\sum_{i=1}^m r_i-\sum_{e\in E(T)}s_e,\\
 e(H)&=\sum_{i=1}^m\left(\binom{r_i}{2}-\binom{|S_i|}{2}\right)
       =\sum_{i=1}^m\binom{r_i}{2}
          -\sum_{e\in E(T)}\binom{s_e}{2}.
 \end{aligned}
\end{equation}
The last equality in each line uses $S_1=\varnothing$ and the
correspondence between each nonroot clique $C_i$ and its edge to
its parent in $T$.

Now fix an integer $x\ge k$ and colour the cliques in the same order.
When $C_i$ is reached, its already coloured vertices are exactly
those in $S_i$, and their colours are distinct. The vertices in $N_i$
must receive distinct colours chosen from the remaining $x-|S_i|$
colours. Thus each colouring of the preceding cliques has
$\prod_{h=|S_i|}^{r_i-1}(x-h)$ extensions to $C_i$. Multiplying
these numbers of choices gives
\[
 P_H(x)=\prod_{i=1}^m\prod_{h=|S_i|}^{r_i-1}(x-h)
       =\prod_{j=1}^k(x-j+1)^{b_j},
\]
where the exponent of $x-j+1$ is
$\#\{i:|S_i|<j\le r_i\}=b_j$ by \eqref{eq:bj}.
Both sides are polynomials in $x$ and agree for every integer
$x\ge k$, so the identity holds for all $x$. In particular, $b_j$
is the multiplicity of $j-1$ as a root of $P_H$. Each $b_j$ is
positive because $H$ contains $K_k$. The colouring argument also
gives $\chi(H)=k$. Since $v(H)=\sum_jb_j$, we obtain
\begin{equation}\label{eq:chordal-target}
 \Phi_H(p)=\prod_{j=1}^k c_j(p)^{b_j},\qquad
 c_j(p)=1-(j-1)(1-p).
\end{equation}

\subsection{An entropy inequality}

We use entropy gluing as developed by Li and Szegedy~\cite{LiSzegedy2011},
Szegedy~\cite{Szegedy2014}, and Lee~\cite{Lee2021}. The construction
below adapts Lee's Markov-tree theorem~\cite[Theorem~2.6]{Lee2021}
to graphons. When all maximal cliques have the same size, the resulting bound
is the clique case of \cite[Theorem~2.7]{Lee2021}. The monotonicity
of $(b_j)$ extends the entropy comparison to cliques of
different sizes.

\begin{proposition}\label{prop:gluing}
Under the hypotheses of \cref{main:profile}, every graphon $W$ satisfies
\begin{equation}\label{eq:general-gluing}
 t(H,W)\prod_{j=1}^{k-1}t(K_j,W)^{b_{j+1}-b_j}
 \ge t(K_k,W)^{b_k}.
\end{equation}
\end{proposition}

\begin{proof}
First assume $\varepsilon\le W\le1$ for some $\varepsilon>0$, so that
all the entropies below are finite.
For a random vector $Z$ with coordinates in $[0,1]$ and probability
density $f$, write
$\mathcal H(Z)=-\int f\log f$ for its entropy. For a graph $F$,
let $\mathrm p_F$ be the probability density on $[0,1]^{V(F)}$
obtained by normalising the product of its edge values. Evaluated
at $Z=(Z_v)_{v\in V(F)}$, it is
\[
 \mathrm p_F\bigl((Z_v)_{v\in V(F)}\bigr)
 =\frac{1}{t(F,W)}\prod_{uv\in E(F)}W(Z_u,Z_v).
\]
Applying Gibbs' inequality to $f$ and $\mathrm p_F$ gives
\begin{equation}\label{eq:entropy-variational}
 \log t(F,W)\ge\mathcal H(Z)
       +\sum_{uv\in E(F)}\E\log W(Z_u,Z_v).
\end{equation}
Equality holds when $Z$ has density $\mathrm p_F$.

Let $(X_1,\ldots,X_k)$ have joint density
\[
 p(x_1, \ldots, x_k) = \frac{1}{t(K_k,W)}\prod_{1\le a<b\le k}W(x_a,x_b).
\]
For $1\le j<k$, its first $j$ coordinates have density
\[
 \frac{1}{t(K_k,W)}
 \int_{[0,1]^{k-j}}\prod_{1\le a<b\le k}W(x_a,x_b)
                 \dd x_{j+1}\cdots\dd x_k.
\]
The joint density is invariant under coordinate permutations, so any
$j$ distinct coordinates have the same law as $(X_1,\ldots,X_j)$.

We will construct a random vector $Y=(Y_v)_{v\in V(H)}$ along the rooted
clique tree. On $C_1$, use the law of $(X_1,\ldots,X_{r_1})$ defined above.
For each subsequent clique $C_i$, order its vertices with those in
$S_i=C_i \cap C_{\pi(i)}$ first. Conditional on their assigned values, sample the
coordinates in $N_i$ according to the conditional law of
$(X_{|S_i|+1},\ldots,X_{r_i})$ given $(X_1,\ldots,X_{|S_i|})$,
\[
 \begin{aligned}
 &p(x_{|S_i|+1},\ldots,x_{r_i}\mid x_1,\ldots,x_{|S_i|})\\
 &\qquad\propto
 \int_{[0,1]^{k-r_i}}\prod_{1\le a<b\le k}W(x_a,x_b)
                 \dd x_{r_i+1}\cdots\dd x_k.
 \end{aligned}
\]
The law already assigned to $S_i$ agrees with this marginal, since
$S_i$ lies in the parent clique. Induction therefore shows that the
coordinates on every $C_i$ have the law of $(X_1,\ldots,X_{r_i})$.
By \eqref{eq:bj}, with $b_{k+1}=0$,
\[
 b_j-b_{j+1}
 =\#\{i:r_i=j\}-\#\{e\in E(T):s_e=j\}.
\]
The entropy chain rule gives
\begin{align*}
 \mathcal H(Y)
 &=\sum_i\mathcal H(X_1,\ldots,X_{r_i})
   -\sum_{e\in E(T)}\mathcal H(X_1,\ldots,X_{s_e})\\
 &=\sum_{j=1}^k\#\{i:r_i=j\}\,\mathcal H(X_1,\ldots,X_j)
   -\sum_{j=1}^k\#\{e\in E(T):s_e=j\}\,\mathcal H(X_1,\ldots,X_j)\\
 &=\sum_{j=1}^k
   \bigl(\#\{i:r_i=j\}-\#\{e\in E(T):s_e=j\}\bigr)
   \mathcal H(X_1,\ldots,X_j)\\
 &=\sum_{j=1}^k(b_j-b_{j+1})\mathcal H(X_1,\ldots,X_j).
\end{align*}
Similarly, \eqref{eq:tree-counts} gives
\begin{align*}
 e(H)
 &=\sum_i\binom{r_i}{2}-\sum_{e\in E(T)}\binom{s_e}{2}\\
 &=\sum_{j=1}^k
   \bigl(\#\{i:r_i=j\}-\#\{e\in E(T):s_e=j\}\bigr)\binom j2\\
 &=\sum_{j=1}^k(b_j-b_{j+1})\binom j2.
\end{align*}
Every edge of $H$ lies in a clique, so its two coordinates in $Y$
have the law of $(X_1,X_2)$. Applying \eqref{eq:entropy-variational}
to $Y$ and substituting these two identities, we obtain
\begin{align}
 \log t(H,W)
 &\ge\mathcal H(Y)+e(H)\E\log W(X_1,X_2)\notag\\
 &=\sum_{j=1}^k(b_j-b_{j+1})
   \left(\mathcal H(X_1,\ldots,X_j)
         +\binom j2\E\log W(X_1,X_2)\right).
 \label{eq:entropy-profile}
\end{align}
Applying \eqref{eq:entropy-variational} to
$(X_1,\ldots,X_j)$ and $K_j$ gives
\[
 \mathcal H(X_1,\ldots,X_j)+\binom j2\E\log W(X_1,X_2)
 \le\log t(K_j,W),
\]
with equality for $j=k$ by the choice of the joint density.
Since $b_j-b_{j+1}\le0$ for $j<k$,
substitution in \eqref{eq:entropy-profile} yields
\[
 \log t(H,W)\ge b_k\log t(K_k,W)
       -\sum_{j=1}^{k-1}(b_{j+1}-b_j)\log t(K_j,W).
\]
Exponentiating proves \eqref{eq:general-gluing}.

For an arbitrary graphon, apply this result to
$W_\varepsilon=\varepsilon+(1-\varepsilon)W$ and let
$\varepsilon\downarrow0$. A telescoping product estimate gives
$|t(F,W_\varepsilon)-t(F,W)|\le e(F)\varepsilon$ for every fixed
finite graph $F$. All powers in \eqref{eq:general-gluing} are
nonnegative integers, so the inequality passes to the limit.
\end{proof}

\subsection{Adjacent clique densities}

The following is the graphon form of the classical clique recurrence
of Moon and Moser~\cite[inequality~(10)]{MoonMoser1962}.

\begin{lemma}\label{lem:mm}
For every $j\ge2$ and every graphon,
\begin{equation}\label{eq:mm}
 jt(K_j,W)^2\le t(K_{j-1},W)t(K_j,W)+(j-1)t(K_{j-1},W)t(K_{j+1},W).
\end{equation}
\end{lemma}

\begin{lemma}\label{lem:clique-step}
If $p>1-1/(k-1)$, then $t(K_j,W)>0$ for $1\le j\le k$ and
\[
 \frac{t(K_j,W)}{t(K_{j-1},W)}\ge 1-(j-1)(1-p),
 \qquad 2\le j\le k.
\]
\end{lemma}

\begin{proof}
We have $t(K_1,W)=1$ and $t(K_2,W)=p>0$, giving the bound for $j=2$.
Suppose $2\le j<k$ and both inequalities hold up to $j$.
Then \eqref{eq:mm} gives
\begin{align*}
 t(K_{j+1},W)
 &\ge\frac{t(K_j,W)}{j-1}
       \left(j\frac{t(K_j,W)}{t(K_{j-1},W)}-1\right)\\
 &\ge\frac{t(K_j,W)}{j-1}
       \bigl(j(1-(j-1)(1-p))-1\bigr)\\
 &=\bigl(1-j(1-p)\bigr)t(K_j,W)>0.
\end{align*}
Here $1-j(1-p)>0$ because $j\le k-1$ and
$p>1-1/(k-1)$. This proves the lemma by induction.
\end{proof}

\begin{proof}[Proof of \cref{main:profile}]
Since $t(K_0,W)=t(K_1,W)=1$, \cref{prop:gluing,lem:clique-step}
give, for $p>1-1/(k-1)$,
\[
 t(H,W)\ge\frac{t(K_k,W)^{b_k}}
                   {\prod_{j<k}t(K_j,W)^{b_{j+1}-b_j}}
       =\prod_{j=1}^k\left(\frac{t(K_j,W)}{t(K_{j-1},W)}\right)^{b_j}
       \ge\prod_{j=1}^k c_j(p)^{b_j}=\Phi_H(p).
\]
\end{proof}

\begin{figure}[htbp]
\centering
\begin{tikzpicture}[
  scale=1,
  vertex/.style={circle,draw=black,fill=white,
                 minimum size=5mm,
                 inner sep=0pt,font=\small},
  clique/.style={rounded corners=6pt,line width=0.6pt,
                 fill opacity=0.12,draw opacity=0.75},
  line width=0.45pt]
\coordinate (v1) at (-1,1.15);
\coordinate (v2) at (-1,-1.15);
\coordinate (v3) at (0,0);
\coordinate (v4) at (1,0);
\coordinate (v5) at (-2,0);
\coordinate (v6) at (2,0.9);
\coordinate (v7) at (2,-0.9);
\path[clique,draw=blue!65!black,fill=blue]
  (-1.4,-1.5) rectangle (1.4,1.5);
\path[clique,draw=orange!80!black,fill=orange]
  (-2.4,-1.7) rectangle (0.4,1.7);
\path[clique,draw=teal!70!black,fill=teal]
  (0.6,-1.3) rectangle (2.4,1.3);
\node[anchor=south,text=blue!65!black] at (0.9,1.5) {$C_2$};
\node[anchor=south,text=orange!80!black] at (-1.9,1.7) {$C_1$};
\node[anchor=south,text=teal!70!black] at (2,1.3) {$C_3$};
\draw
      (v1)--(v2) (v1)--(v3) (v1)--(v4)
      (v2)--(v3) (v2)--(v4) (v3)--(v4)
      (v1)--(v5) (v2)--(v5) (v3)--(v5)
      (v4)--(v6) (v4)--(v7) (v6)--(v7);
\foreach \v in {1,...,7}
  \node[vertex] at (v\v) {$\v$};
\end{tikzpicture}
\caption{The chordal graph $H$ and its maximal cliques $C_1,C_2,C_3$.}
\label{fig:chordal-example}
\end{figure}
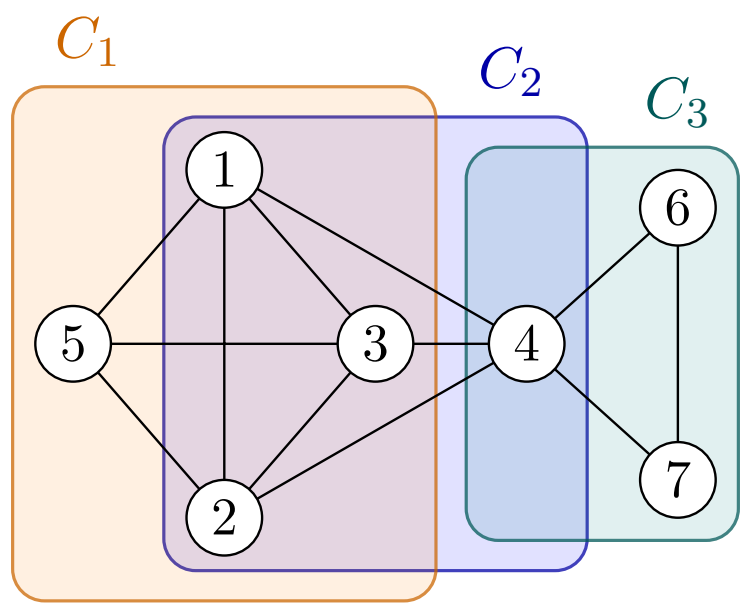

The graph $H$ in \cref{fig:chordal-example} has maximal cliques
$C_1=\{1,2,3,5\}$, $C_2=\{1,2,3,4\}$, and $C_3=\{4,6,7\}$.
The path $C_1$-$C_2$-$C_3$ is a clique tree, with clique sizes $4,4,3$
and intersections of sizes $3,1$. Thus
$(b_1,b_2,b_3,b_4)=(1,2,2,2)$, and \cref{main:profile} proves the
Goodman-type lower bound for $H$.

\section{Universal vertices}\label{sec:lifting}

We prove that adding a universal vertex preserves the bound when
$\Phi_F$ is nonnegative, nondecreasing, and convex on the relevant
interval. The proof begins with a weighted inequality for triangle
densities.

\subsection{A weighted triangle inequality}

For a graphon $W$, let $d_W(x)=\int W(x,y)\dd y$ be its degree
function. Set
\begin{align*}
 \tau_W(x)&=\int W(x,y)W(x,z)W(y,z)\dd y\dd z,\\
 M_j(W)&=\int d_W(x)^j\dd x\quad(j\ge0).
\end{align*}
Thus $\tau_W(x)$ is the triangle density in $W$ with one vertex fixed at $x$.
For a graphon $V$, let $(T_Vf)(x)=\int V(x,y)f(y)\dd y$ denote
its integral operator. We use the inner product
$\inner{f}{g}=\int f(x)g(x)\dd x$ on real $L^2[0,1]$.

\begin{lemma}\label{lem:weighted-triangle}
If $p=\int W(x,y)\dd x\dd y>0$, then, for every integer $h\ge0$,
\begin{equation}\label{eq:weighted-triangle}
 \int d_W(x)^h\tau_W(x)\dd x\ge\left(2-\frac1p\right)M_{h+2}(W).
\end{equation}
\end{lemma}

\begin{proof}
Applying $ab\ge a+b-1$ to $a=W(x,y)$ and $b=W(x,z)$, and then
using the symmetry of $W$, we obtain
\begin{align*}
 \tau_W(x)
 &=\int W(x,y)W(x,z)W(y,z)\dd y\dd z\\
 &\ge\int\bigl(W(x,y)+W(x,z)-1\bigr)W(y,z)\dd y\dd z\\
 &=\int W(x,y)d_W(y)\dd y+\int W(x,z)d_W(z)\dd z-p\\
 &=2(T_Wd_W)(x)-p.
\end{align*}
Multiplying by $d_W(x)^h\ge0$ and integrating in $x$ gives
\begin{align*}
 \int d_W(x)^h\tau_W(x)\dd x
 &\ge2\int d_W(x)^h(T_Wd_W)(x)\dd x
       -p\int d_W(x)^h\dd x\\
 &=2\inner{d_W^h}{T_Wd_W}-pM_h(W).
\end{align*}
Set $U=1-W$ and $q=1-p$, so that $d_U=1-d_W$. Since
$T_Wd_W=d_W-q+T_Ud_U$,
\begin{align*}
 \inner{d_W^h}{T_Wd_W}
 &=\int d_W(x)^h(T_Wd_W)(x)\dd x\\
 &=\int d_W(x)^h\bigl(d_W(x)-q+(T_Ud_U)(x)\bigr)\dd x\\
 &=\int d_W(x)^{h+1}\dd x-q\int d_W(x)^h\dd x
       +\int d_W(x)^h(T_Ud_U)(x)\dd x\\
 &=M_{h+1}(W)-qM_h(W)+\inner{d_W^h}{T_Ud_U}.
\end{align*}
Since $d_U(x)=\int U(x,y)\dd y$, symmetry of $U$ gives
\begin{align*}
 &2\left(\inner{d_W^h}{T_Ud_U}-\int d_W(x)^h d_U(x)^2\dd x\right)\\
 &\qquad=2\int U(x,y)d_W(x)^h d_U(y)\dd x\dd y
          -2\int U(x,y)d_W(x)^h d_U(x)\dd x\dd y\\
 &\qquad=2\int U(x,y)d_W(x)^h(d_U(y)-d_U(x))\dd x\dd y\\
 &\qquad=\int U(x,y)d_W(x)^h(d_U(y)-d_U(x))\dd x\dd y\\
 &\qquad\quad+\int U(x,y)d_W(y)^h(d_U(x)-d_U(y))\dd x\dd y\\
 &\qquad=\int U(x,y)(d_U(y)-d_U(x))
                 (d_W(x)^h-d_W(y)^h)\dd x\dd y\\
 &\qquad\ge0,
\end{align*}
where the last inequality follows from
$d_U(y)-d_U(x)=d_W(x)-d_W(y)$ and the monotonicity of $t\mapsto t^h$
on $[0,1]$. Since $d_U=1-d_W$,
\begin{align*}
 \int d_W(x)^h d_U(x)^2\dd x
 &=\int d_W(x)^h(1-d_W(x))^2\dd x\\
 &=M_h(W)-2M_{h+1}(W)+M_{h+2}(W).
\end{align*}
Combining these estimates and using $q=1-p$, we obtain
\begin{align*}
 &\int d_W(x)^h\tau_W(x)\dd x\\
 &\quad\ge2\inner{d_W^h}{T_Wd_W}-pM_h(W)\\
 &\quad=2M_{h+1}(W)-(2q+p)M_h(W)+2\inner{d_W^h}{T_Ud_U}\\
 &\quad\ge2M_{h+1}(W)-(2q+p)M_h(W)
           +2\int d_W(x)^h d_U(x)^2\dd x\\
 &\quad=2M_{h+1}(W)-(2q+p)M_h(W)
           +2\bigl(M_h(W)-2M_{h+1}(W)+M_{h+2}(W)\bigr)\\
 &\quad=2M_{h+2}(W)-2M_{h+1}(W)+pM_h(W)\\
 &\quad=\left(2-\frac1p\right)M_{h+2}(W)
           +\frac1p\bigl(M_{h+2}(W)-2pM_{h+1}(W)+p^2M_h(W)\bigr)\\
 &\quad=\left(2-\frac1p\right)M_{h+2}(W)
           +\frac1p\int d_W(x)^h(d_W(x)-p)^2\dd x\\
 &\quad\ge\left(2-\frac1p\right)M_{h+2}(W).
\end{align*}
This proves the claim.
\end{proof}

\subsection{Adding a universal vertex}

\begin{proposition}\label{prop:lifting}
Let $F$ have $n\ge2$ vertices, and let $0\le\underline{p}<1$.
Suppose $\Phi_F$ is nondecreasing and convex on $[\underline{p},1]$, with
$\Phi_F(\underline{p})=0$, and $t(F,V)\ge\Phi_F(z)$ for every graphon $V$
of edge density $z>\underline{p}$.
Then every graphon $W$ of edge density $p>1/(2-\underline{p})$ satisfies
\[
 t(K_1\join F,W)\ge\Phi_{K_1\join F}(p).
\]
Moreover, $\Phi_{K_1\join F}$ is nondecreasing and convex on
$[1/(2-\underline{p}),1]$, with $\Phi_{K_1\join F}(1/(2-\underline{p}))=0$.
Thus $K_1\join F$ satisfies the hypotheses of the proposition with
$1/(2-\underline{p})$ in place of $\underline{p}$.
\end{proposition}

\begin{proof}
For any $c\in[\underline{p},1]$, the tangent line to $\Phi_F$ at $c$ gives
\begin{equation}\label{eq:tangent}
 t(F,V)\ge\Phi_F(c)+\Phi_F'(c)(z-c),\qquad 0\le z\le1.
\end{equation}
For $z>\underline{p}$, this follows from the assumed bound and convexity.
For $z\le\underline{p}$, nonnegativity of $\Phi_F'(c)$ and convexity give
\begin{align*}
 \Phi_F(c)+\Phi_F'(c)(z-c)
 &\le\Phi_F(c)+\Phi_F'(c)(\underline{p}-c)\\
 &\le\Phi_F(\underline{p})=0\le t(F,V).
\end{align*}

By definition,
\[
 t(K_1\join F,W)
 =\int\left(\int\prod_{uv\in E(F)}W(y_u,y_v)
      \prod_{v\in V(F)}W(x,y_v)\prod_{v\in V(F)}\dd y_v\right)\dd x.
\]
Fix $x\in[0,1]$ with $d_W(x)>0$, and define the map $Q_x:[0,1]\to[0,1]$ by
\[
 Q_x(s):=\inf\left\{y\in[0,1]:
       \int_0^y W(x,z)\dd z\ge s\,d_W(x)\right\},\qquad s\in[0,1].
\]
The map $Q_x$ is nondecreasing and satisfies
\[
 \int\ind_{\{Q_x(s)\le y\}}\dd s
 =\frac{1}{d_W(x)}\int_0^y W(x,z)\dd z,
 \qquad y\in[0,1].
\]
Consequently, for every bounded measurable $f:[0,1]\to\mathbb R$,
\[
 \int f(Q_x(s))\dd s
 =\frac{1}{d_W(x)}\int f(y)W(x,y)\dd y.
\]
Define the reweighted graphon
\[
 W_x(s,t)=W(Q_x(s),Q_x(t)).
\]
Applying the integration identity for $Q_x$ successively to the variables
$(s_v)_{v\in V(F)}$, we obtain
\begin{align*}
 t(F,W_x)
 &=\int\prod_{uv\in E(F)}W_x(s_u,s_v)\prod_{v\in V(F)}\dd s_v\\
 &=\int\prod_{uv\in E(F)}W(Q_x(s_u),Q_x(s_v))
       \prod_{v\in V(F)}\dd s_v\\
 &=\int\prod_{uv\in E(F)}W(y_u,y_v)
       \left(\prod_{v\in V(F)}\frac{W(x,y_v)}{d_W(x)}\right)
       \prod_{v\in V(F)}\dd y_v\\
 &=\frac{1}{d_W(x)^n}\int\prod_{uv\in E(F)}W(y_u,y_v)
       \prod_{v\in V(F)}W(x,y_v)\prod_{v\in V(F)}\dd y_v.
\end{align*}
In particular, its edge density is
\[
 t(K_2,W_x)
 =\frac{1}{d_W(x)^2}\int W(y,z)W(x,y)W(x,z)\dd y\dd z
 =\frac{\tau_W(x)}{d_W(x)^2}.
\]
Multiplying the expression for $t(F,W_x)$ by $d_W(x)^n$ and integrating in $x$ gives
\[
 t(K_1\join F,W)=\int d_W(x)^n t(F,W_x)\dd x.
\]
When $d_W(x)=0$, the contribution of $x$ is zero; we set $W_x=0$ there.

\Cref{lem:weighted-triangle}, with $h=n-2$, gives
\begin{align*}
 \int d_W(x)^n t(K_2,W_x)\dd x
 &=\int d_W(x)^{n-2}\tau_W(x)\dd x\\
 &\ge\left(2-\frac1p\right)M_n(W).
\end{align*}
For $p>1/(2-\underline{p})$, the point $2-1/p$ lies in $(\underline{p},1]$.
Applying \eqref{eq:tangent} at this point to each $W_x$ and integrating gives
\begin{align*}
 &t(K_1\join F,W)\\
 &\quad=\int d_W(x)^n t(F,W_x)\dd x\\
 &\quad\ge\int d_W(x)^n
     \left[\Phi_F\!\left(2-\frac1p\right)
       +\Phi_F'\!\left(2-\frac1p\right)
          \left(t(K_2,W_x)-\left(2-\frac1p\right)\right)\right]\dd x\\
 &\quad=M_n(W)\Phi_F\!\left(2-\frac1p\right)
       +\Phi_F'\!\left(2-\frac1p\right)
          \left(\int d_W(x)^n t(K_2,W_x)\dd x
                     -\left(2-\frac1p\right)M_n(W)\right)\\
 &\quad\ge M_n(W)\Phi_F\!\left(2-\frac1p\right)\\
 &\quad\ge\left(\int d_W(x)\dd x\right)^n
              \Phi_F\!\left(2-\frac1p\right)\\
 &\quad=p^n\Phi_F\!\left(2-\frac1p\right).
\end{align*}
Here we use $\Phi_F'(2-1/p),\Phi_F(2-1/p)\ge0$ and the preceding bound
on $\int d_W(x)^n t(K_2,W_x)\dd x$; the last inequality follows from
Jensen's inequality.
Since $P_{K_1\join F}(x)=xP_F(x-1)$,
\[
 \Phi_{K_1\join F}(p)=p^n\Phi_F\!\left(2-\frac1p\right).
\]
Thus $K_1\join F$ also satisfies the lower bound for $p>1/(2-\underline{p})$.

For the second assertion, differentiating gives
\begin{align*}
 \Phi_{K_1\join F}'(p)
 &=np^{n-1}\Phi_F\!\left(2-\frac1p\right)
   +p^{n-2}\Phi_F'\!\left(2-\frac1p\right),\\
 \Phi_{K_1\join F}''(p)
 &=n(n-1)p^{n-2}\Phi_F\!\left(2-\frac1p\right)
   +2(n-1)p^{n-3}\Phi_F'\!\left(2-\frac1p\right)\\
 &\qquad+p^{n-4}\Phi_F''\!\left(2-\frac1p\right).
\end{align*}
These are nonnegative on $[1/(2-\underline{p}),1]$, since $2-1/p$ lies
in $[\underline{p},1]$ there and $\Phi_F,\Phi_F',\Phi_F''\ge0$ on
$[\underline{p},1]$. Finally,
\[
 \Phi_{K_1\join F}\!\left(\frac1{2-\underline{p}}\right)
 =\frac{\Phi_F(\underline{p})}{(2-\underline{p})^n}=0.
\]
\end{proof}

Iterating \cref{prop:lifting} shows that, under its hypotheses,
$t(K_s\join F,W)\ge\Phi_{K_s\join F}(p)$ for every integer $s\ge0$
and every graphon $W$ of edge density
\begin{equation}\label{eq:iterated-threshold}
 p>1-\frac1{s+1/(1-\underline{p})}.
\end{equation}
Indeed, with $\underline{p}_0=\underline{p}$, the thresholds obey
$\underline{p}_{j+1}=1/(2-\underline{p}_j)$, or equivalently
$(1-\underline{p}_{j+1})^{-1}=(1-\underline{p}_j)^{-1}+1$;
by the second part of \cref{prop:lifting}, the hypotheses are
preserved at every step.

\begin{proof}[Proof of \cref{main:odd-join}]
For odd $\ell\ge3$,
\[
 \Phi_{C_\ell}(p)=p\bigl(p^{\ell-1}-(1-p)^{\ell-1}\bigr).
\]
This polynomial is nonnegative on $[1/2,1]$, with $\Phi_{C_\ell}(1/2)=0$.
On the same interval, its derivatives satisfy
\begin{align*}
 \Phi_{C_\ell}'(p)
 &=p^{\ell-1}-(1-p)^{\ell-1}
   +(\ell-1)p\bigl(p^{\ell-2}+(1-p)^{\ell-2}\bigr)\ge0,\\
 \Phi_{C_\ell}''(p)
 &=2(\ell-1)\bigl(p^{\ell-2}+(1-p)^{\ell-2}\bigr)\\
 &\qquad+(\ell-1)(\ell-2)p\bigl(p^{\ell-3}-(1-p)^{\ell-3}\bigr)\ge0.
\end{align*}
Apply \cref{main:odd,prop:lifting} and
\eqref{eq:iterated-threshold} with $\underline{p}=1/2$.
Since $\chi(K_s\join C_\ell)=s+3$, the resulting threshold is
the required one.
\end{proof}

\subsection{Trivially perfect graphs}\label{sec:tp}

\begin{proof}[Proof of \cref{main:tp}]
We use strong induction on $v(H)$, assuming that every trivially perfect
graph $F$ with at least one edge and $v(F)<v(H)$ satisfies the
Goodman-type lower bound, with $\Phi_F$ nondecreasing and convex on
$[\underline{p}_F,1]$ and $\Phi_F(\underline{p}_F)=0$.

The base case is a star $H$ with $m\ge1$ leaves. Here $\underline{p}_H=0$,
and Jensen's inequality gives
\[
 t(H,W)=\int d_W(x)^m\dd x\ge p^m=\Phi_H(p).
\]
The polynomial $p^m$ is nondecreasing and convex on $[0,1]$ and equals zero at $p=0$.

If $H$ is disconnected, write $H=H_1\sqcup H_2$, where $H_1,H_2$
are nonempty. Since $\chi(H)=\max\{\chi(H_1),\chi(H_2)\}$,
\[
 t(H,W)=t(H_1,W)t(H_2,W)
 \ge\Phi_{H_1}(p)\Phi_{H_2}(p)=\Phi_H(p),\qquad p>\underline{p}_H.
\]
Here $P_H=P_{H_1}P_{H_2}$ and $v(H)=v(H_1)+v(H_2)$ give the last
equality; if $H_i$ is edgeless, then $t(H_i,W)=\Phi_{H_i}(p)=1$.
The product rule gives nonnegativity of $\Phi_H'$ and $\Phi_H''$ on
$[\underline{p}_H,1]$. Relabelling so that $\chi(H_1)=\chi(H)$, we have
\[
 \Phi_H(\underline{p}_H)
 =\Phi_{H_1}(\underline{p}_{H_1})\Phi_{H_2}(\underline{p}_H)=0.
\]

It remains to consider $H=K_1\join F$, where $F$ is trivially perfect
and has an edge~\cite{DrangePilipczuk}. The induction hypothesis allows us to apply
\cref{prop:lifting} with $\underline{p}=\underline{p}_F$.
This gives the required bound and the required properties of
$\Phi_H$, since $\chi(H)=\chi(F)+1$ and
\[
 \frac1{2-\underline{p}_F}=1-\frac1{\chi(F)}=\underline{p}_H.
\]
\end{proof}

\section{Eventual validity and girth}\label{sec:eventual}

We first expand $t(H,1-U)-\Phi_H(1-q)$ over forests, where
$q=t(K_2,U)$. When the girth is odd, estimates for tree and cycle
densities give the lower bound for small maximum degree. We then
pass from small maximum degree to small edge density. The section
ends with the obstruction for even girth and the proof of
\cref{main:eventual}.

For a graphon $U$ of edge density $q$, we use
\begin{equation}\label{evt:eq:degree-notation}
 \begin{aligned}
 d_U(x)&=\int U(x,y)\dd y,
 &\Delta_U&=\norm{d_U}_\infty,\\
 M_j(U)&=\int d_U(x)^j\dd x\quad(j\ge0),
 &\operatorname{Var}(d_U)&=\int(d_U(x)-q)^2\dd x=M_2(U)-q^2.
 \end{aligned}
\end{equation}
Here $\operatorname{Var}(d_U)$ is the variance of the degree function,
and the norm is the essential supremum. Thus $0\le q\le\Delta_U\le1$.

\subsection{Expansion over forests}\label{evt:sec:expansion}

For $A\subseteq E(H)$, let $H_A=(V(H),A)$ and define
\[
 r(A)=v(H)-\comp(H_A),
\]
where $\comp(H_A)$ counts all connected components, including isolated
vertices. Thus $r(A)$ is the number of edges in a spanning forest of $H_A$.

\begin{lemma}\label{evt:lem:expansion}
For every graph $H$ and every graphon $U$ of edge density $q$,
\begin{align}
 \Phi_H(1-q)&=\sum_{A\subseteq E(H)}(-1)^{|A|}q^{r(A)},
 \label{evt:eq:chromatic-expansion}\\
 t(H,1-U)&=\sum_{A\subseteq E(H)}(-1)^{|A|}t(H_A,U).
 \label{evt:eq:complement-expansion}
\end{align}
\end{lemma}

\begin{proof}
Whitney's spanning-subgraph expansion~\cite{Whitney1931} follows by
inclusion--exclusion over monochromatic edges:
\[
 P_H(k)=\sum_{A\subseteq E(H)}(-1)^{|A|}k^{\comp(H_A)}
 \qquad(k\in\mathbb N).
\]
Both sides are polynomials in $k$, so the identity holds for all $k \in \mathbb{R}$.
For $q>0$, substituting $k=1/q$ gives
\[
 \Phi_H(1-q)
 =q^{v(H)}P_H(1/q)
 =\sum_{A\subseteq E(H)}(-1)^{|A|}q^{v(H)-\comp(H_A)}.
\]
Continuity gives the same formula at $q=0$. Also,
\begin{align*}
 t(H,1-U)
 &=\int\prod_{uv\in E(H)}(1-U(x_u,x_v))\prod_{v\in V(H)}\dd x_v\\
 &=\sum_{A\subseteq E(H)}(-1)^{|A|}
       \int\prod_{uv\in A}U(x_u,x_v)\prod_{v\in V(H)}\dd x_v\\
 &=\sum_{A\subseteq E(H)}(-1)^{|A|}t(H_A,U).
\end{align*}
\end{proof}

Subtracting the two identities gives
\[
 t(H,1-U)-\Phi_H(1-q)
 =\sum_{A\subseteq E(H)}(-1)^{|A|}
       \bigl(t(H_A,U)-q^{r(A)}\bigr).
\]
The summands with $r(A)=0$ or $r(A)=1$ are zero: the empty edge set
contributes $1-1=0$, and each single edge contributes $-(q-q)=0$.

Fix a total order $<$ on $E(H)$, and regard each forest
$F\subseteq E(H)$ as the spanning graph $(V(H),F)$. Let $B(F)$ be
the set of edges $e\in E(H)\setminus F$ such that adding $e$ to $F$
creates a cycle and every other edge of this cycle precedes $e$
in the fixed order. Set
\[
 I_F(U)=\int\prod_{uv\in F}U(x_u,x_v)
       \prod_{uv\in B(F)}(1-U(x_u,x_v))\prod_{v\in V(H)}\dd x_v.
\]

\begin{lemma}\label{evt:lem:forest-expansion}
For every graphon $U$ of edge density $q$,
\begin{align}
 \Phi_H(1-q)&=\sum_{\substack{F\text{ forest}\\B(F)=\varnothing}}
                   (-1)^{|F|}q^{|F|},
 \label{evt:eq:forest-chromatic}\\
 t(H,1-U)&=\sum_{F\text{ forest}}(-1)^{|F|}I_F(U).
 \label{evt:eq:forest-density}
\end{align}
Here the sums range over forests $F\subseteq E(H)$.
\end{lemma}

Identity \eqref{evt:eq:forest-chromatic} is also known as Whitney's broken-circuit
theorem~\cite{Whitney1932}, where a \emph{broken circuit} is a cycle with
its largest edge removed; $B(F)=\varnothing$ if and only if $F$
contains no broken circuit.

\begin{proof}
We first show that every $A\subseteq E(H)$ has a unique representation
\[
 A=F\mathbin{\dot\cup}D,\qquad F\text{ a forest},\quad D\subseteq B(F).
\]

Assign the weights $2,4,\ldots,2^{e(H)}$ to the edges of $H$ in the fixed
order. Given $A\subseteq E(H)$, let $F$ be the minimum spanning
forest of $H_A$. It is unique because different edge sets have
different total weights. Adding an edge $e\in A\setminus F$ to $F$
creates a cycle. Every other edge in this cycle has smaller weight
than $e$: otherwise, replacing a heavier edge by $e$ would give a
spanning forest of smaller total weight. Hence $e\in B(F)$, giving
the required representation with $D=A\setminus F$.

For uniqueness, suppose $A=F\mathbin{\dot\cup}D$ with $F$ a forest and
$D\subseteq B(F)$. Since every edge of $D$ has its endpoints in the
same component of $F$, the graphs $F$ and $H_A$ have the same
components. Let $F'$ be the minimum spanning forest of $H_A$.
If $F\ne F'$, choose $e\in F'\setminus F\subseteq D$.
Adding $e$ to $F$ creates a cycle whose other edges all have smaller
weight than $e$. As the weights are distinct powers of $2$, the
sum of their weights is smaller than the weight of $e$.
Replacing $e$ in $F'$ by these other edges gives a graph with the
same components as $H_A$ and smaller total weight than $F'$.
A spanning forest of this graph would also have smaller total
weight than $F'$, contradicting its minimality. Thus $F=F'$, so
$A$ determines $F$ and then $D=A\setminus F$. We also have
$r(A)=|F|$.

Using this decomposition in the expansion of $t(H,1-U)$, we obtain
\begin{align*}
 &t(H,1-U)\\
 &\quad=\sum_{F\text{ forest}}\sum_{D\subseteq B(F)}(-1)^{|F|+|D|}
       \int\prod_{uv\in F\mathbin{\dot\cup}D}U(x_u,x_v)\prod_{v\in V(H)}\dd x_v\\
 &\quad=\sum_{F\text{ forest}}(-1)^{|F|}
       \int\prod_{uv\in F}U(x_u,x_v)
       \prod_{uv\in B(F)}(1-U(x_u,x_v))\prod_{v\in V(H)}\dd x_v\\
 &\quad=\sum_{F\text{ forest}}(-1)^{|F|}I_F(U).
\end{align*}
For each forest $F$, the corresponding terms in
\eqref{evt:eq:chromatic-expansion} sum to
\[
 (-1)^{|F|}q^{|F|}\sum_{D\subseteq B(F)}(-1)^{|D|}
 =\begin{cases}
   (-1)^{|F|}q^{|F|},&B(F)=\varnothing,\\
   0,&B(F)\ne\varnothing.
  \end{cases}
\]
Summing over $F$ proves \eqref{evt:eq:forest-chromatic}.
\end{proof}

Since $I_F(U)=t(F,U)$ when $B(F)=\varnothing$, subtracting the two
identities gives
\begin{align}
 t(H,1-U)-\Phi_H(1-q)
 &=\sum_{\substack{F\text{ forest}\\B(F)=\varnothing}}
       (-1)^{|F|}\bigl(t(F,U)-q^{|F|}\bigr)\notag\\
 &\qquad+\sum_{\substack{F\text{ forest}\\B(F)\ne\varnothing}}
       (-1)^{|F|}I_F(U).
 \label{evt:eq:forest-gamma}
\end{align}

\subsection{Tree and forest densities}\label{evt:sec:forests}

The lower bound below is the tree case of Sidorenko's
inequality~\cite{Sidorenko1993}; the upper bound is his comparison
with the star having the same number of edges~\cite{Sidorenko1994}.

\begin{lemma}\label{evt:lem:tree}
If $T$ is a tree with $r\ge1$ edges, then
\begin{equation}\label{evt:eq:tree-moments}
 q^r\le t(T,U)\le M_r(U).
\end{equation}
\end{lemma}

Sidorenko's inequality for trees holds with equality whenever $U$ is
regular, that is, $d_U(x)=q$ almost everywhere. The following lemma
bounds the gap $t(F,U)-q^{e(F)}$ for forests $F$ in terms of the
variance of the degree function.

\begin{lemma}\label{evt:lem:forest}
If $F$ is a forest with $r\ge2$ edges, then
\begin{equation}\label{evt:eq:forest-variance}
 0\le t(F,U)-q^r\le\binom r2\Delta_U^{r-2}\operatorname{Var}(d_U).
\end{equation}
In particular, taking $F$ to be the star $K_{1,r}$ gives
\begin{equation}\label{evt:eq:moment-variance}
 0\le M_r(U)-q^r\le\binom r2\Delta_U^{r-2}\operatorname{Var}(d_U).
\end{equation}
\end{lemma}

\begin{proof}
Assume $\Delta_U>0$, since the case $\Delta_U=0$ is immediate.
Let $T_1,\ldots,T_\ell$ be the components of $F$ containing an edge,
and write $r_i=e(T_i)$, so that $\sum_i r_i=r$.
By \cref{evt:lem:tree} and H\"older's inequality,
\begin{align*}
 q^r=\prod_{i=1}^{\ell}q^{r_i}
 &\le t(F,U)=\prod_{i=1}^{\ell}t(T_i,U)\\
 &\le\prod_{i=1}^{\ell}M_{r_i}(U)
 \le\prod_{i=1}^{\ell}M_r(U)^{r_i/r}=M_r(U).
\end{align*}
On $[0,\Delta_U]$, the function $s\mapsto s^r$ has second
derivative at most $r(r-1)\Delta_U^{r-2}$.
Since $q,d_U(x)\in[0,\Delta_U]$ for almost every $x$, Taylor's
theorem at $s=q$ gives
\[
 d_U(x)^r\le q^r+rq^{r-1}(d_U(x)-q)
       +\binom r2\Delta_U^{r-2}(d_U(x)-q)^2.
\]
Integrating and using $\int(d_U(x)-q)\dd x=0$, we obtain
\[
 0\le t(F,U)-q^r\le M_r(U)-q^r
 \le\binom r2\Delta_U^{r-2}\operatorname{Var}(d_U).
\]
The star case follows from $t(K_{1,r},U)=M_r(U)$.
\end{proof}

\subsection{Path and cycle densities}\label{evt:sec:triangles}

We bound the difference between the densities of a path and the
cycle obtained by joining its endpoints.

\begin{lemma}\label{evt:lem:spanning}
Let $F$ be a forest with $r$ edges, and let $u,v\in V(F)$ have
distance $s-1\ge2$. Then
\begin{equation}\label{evt:eq:spanning}
 \begin{aligned}
 0&\le\int\biggl(\prod_{ab\in E(F)}U(x_a,x_b)\biggr)
       (1-U(x_u,x_v))\prod_{w\in V(F)}\dd x_w\\
 &\le\Delta_U^{r-s+1}\bigl(t(P_s,U)-t(C_s,U)\bigr).
 \end{aligned}
\end{equation}
\end{lemma}

\begin{proof}
Write the $u$--$v$ path as $u=v_1,\ldots,v_s=v$.
If $a$ is a leaf outside this path with neighbour $b$, then $x_a$
occurs only in the factor $U(x_a,x_b)$. Hence
\begin{align*}
 &\int\biggl(\prod_{cd\in E(F)}U(x_c,x_d)\biggr)(1-U(x_u,x_v))\dd x_a\\
 &\quad=d_U(x_b)\biggl(\prod_{cd\in E(F)\setminus\{ab\}}U(x_c,x_d)\biggr)
       (1-U(x_u,x_v))\\
 &\quad\le\Delta_U\biggl(\prod_{cd\in E(F)\setminus\{ab\}}U(x_c,x_d)\biggr)
       (1-U(x_u,x_v)).
\end{align*}
Repeating this for the $r-s+1$ edges outside the path and
integrating isolated vertices, which contribute $1$, gives
\begin{align*}
 0&\le\int\biggl(\prod_{ab\in E(F)}U(x_a,x_b)\biggr)
       (1-U(x_u,x_v))\prod_{w\in V(F)}\dd x_w\\
 &\le\Delta_U^{r-s+1}
       \int\biggl(\prod_{i=1}^{s-1}U(x_{v_i},x_{v_{i+1}})\biggr)
       (1-U(x_u,x_v))\prod_{i=1}^s\dd x_{v_i}\\
 &=\Delta_U^{r-s+1}\bigl(t(P_s,U)-t(C_s,U)\bigr).
\end{align*}
\end{proof}

For every $s\ge3$,
\[
 t(P_s,U)-t(C_s,U)
 =\int\biggl(\prod_{i=1}^{s-1}U(x_i,x_{i+1})\biggr)
       (1-U(x_1,x_s))\prod_{i=1}^s\dd x_i.
\]
For $s=3$,
\begin{equation}\label{evt:eq:Y}
 \begin{aligned}
 M_2(U)-t(C_3,U)
 &=t(P_3,U)-t(C_3,U)\\
 &=\int U(x,y)U(y,z)(1-U(x,z))\dd x\dd y\dd z.
 \end{aligned}
\end{equation}

\begin{lemma}\label{evt:lem:triangle-path}
For every integer $s\ge3$,
\begin{equation}\label{evt:eq:triangle-path}
 0\le t(P_s,U)-t(C_s,U)
 \le(s-2)\Delta_U^{s-3}\bigl(M_2(U)-t(C_3,U)\bigr).
\end{equation}
\end{lemma}

\begin{proof}
For $s=3$, the claim is immediate by \eqref{evt:eq:Y}.

Define the iterated kernels recursively by
\[
 U^{\circ1}(x,y)=U(x,y),\qquad
 U^{\circ(\ell+1)}(x,y)=\int U^{\circ\ell}(x,z)U(z,y)\dd z.
\]
Induction gives
\begin{equation}\label{evt:eq:path-kernel}
 0\le U^{\circ\ell}(x,y)\le\Delta_U^{\ell-1}
 \quad\text{for almost every }(x,y),
\end{equation}
since $U\le1$ and
\[
 U^{\circ(\ell+1)}(x,y)
 \le\Delta_U^{\ell-1}\int U(z,y)\dd z
 =\Delta_U^{\ell-1}d_U(y)\le\Delta_U^\ell.
\]

We obtain
\begin{align*}
 &t(P_s,U)-t(C_s,U)\\
 &\quad=\int\biggl(\prod_{i=1}^{s-1}U(x_i,x_{i+1})\biggr)
       (1-U(x_1,x_s))\prod_{i=1}^s\dd x_i\\
 &\quad\le\int\biggl(\prod_{i=1}^{s-1}U(x_i,x_{i+1})\biggr)
       (1-U(x_1,x_{s-1}))\prod_{i=1}^s\dd x_i\\
 &\qquad\quad+\int\biggl(\prod_{i=1}^{s-1}U(x_i,x_{i+1})\biggr)
       U(x_1,x_{s-1})(1-U(x_1,x_s))\prod_{i=1}^s\dd x_i\\
 &\quad=\int\biggl(\prod_{i=1}^{s-2}U(x_i,x_{i+1})\biggr)
       d_U(x_{s-1})(1-U(x_1,x_{s-1}))\prod_{i=1}^{s-1}\dd x_i\\
 &\qquad\quad+\int U^{\circ(s-2)}(x,y)U(x,y)U(y,z)
       (1-U(x,z))\dd x\dd y\dd z\\
 &\quad\le\Delta_U\bigl(t(P_{s-1},U)-t(C_{s-1},U)\bigr)
       +\Delta_U^{s-3}\bigl(M_2(U)-t(C_3,U)\bigr),
\end{align*}
where the first inequality follows from
\[
 1-U(x_1,x_s)
 \le1-U(x_1,x_{s-1})
       +U(x_1,x_{s-1})(1-U(x_1,x_s)),
\]
and the last follows from $d_U\le\Delta_U$,
\eqref{evt:eq:path-kernel}, and \eqref{evt:eq:Y}.
In the second equality, we integrate $x_s$ in the first term and
$x_2,\ldots,x_{s-2}$ in the second.
Induction on $s$ proves \eqref{evt:eq:triangle-path}.
\end{proof}

\subsection{Comparison of cycle densities}\label{evt:sec:spectral}

For odd $g\ge5$, the following comparison provides the estimate
needed for cycles of length at least $g$.

\begin{theorem}\label{evt:thm:cycle-comparison}
Let $g\ge5$ be odd and let $s\ge g$ be an integer. For every graphon $U$,
\begin{equation}\label{evt:eq:cycle-comparison}
 \begin{aligned}
 0&\le M_{s-1}(U)-t(C_s,U)\\
 &\le\frac{(s-2)(s-3)}{(g-2)(g-3)}\,
       \Delta_U^{s-g}\bigl(M_{g-1}(U)-t(C_g,U)\bigr).
 \end{aligned}
\end{equation}
\end{theorem}

\begin{lemma}\label{evt:lem:scalar}
Let $2\le a\le b$ be integers and let $\Delta>0$. If a random
variable $Z$ takes values in $[0,\Delta]$ and $0\le\alpha\le\E Z$, then
\begin{equation}\label{evt:eq:scalar}
 \E Z^b-\alpha^b
 \le\frac{b(b-1)}{a(a-1)}\Delta^{b-a}
       (\E Z^a-\alpha^a).
\end{equation}
\end{lemma}

\begin{proof}
Equality holds when $a=b$. For $b>a$, consider the polynomial
\[
 z^b-\frac{b(b-1)}{a(a-1)}\Delta^{b-a}z^a.
\]
Its first and second derivatives on $[0,\Delta]$ are, respectively,
\begin{align*}
 bz^{a-1}\left(z^{b-a}-\frac{b-1}{a-1}\Delta^{b-a}\right)&\le0,\\
 b(b-1)z^{a-2}\bigl(z^{b-a}-\Delta^{b-a}\bigr)&\le0.
\end{align*}
The polynomial is therefore nonincreasing and concave. Jensen's
inequality and $\alpha\le\E Z$ give
\begin{align*}
 \E Z^b-\frac{b(b-1)}{a(a-1)}\Delta^{b-a}\E Z^a
 &\le(\E Z)^b-\frac{b(b-1)}{a(a-1)}\Delta^{b-a}(\E Z)^a\\
 &\le\alpha^b-\frac{b(b-1)}{a(a-1)}\Delta^{b-a}\alpha^a.
\end{align*}
Rearranging proves the claim.
\end{proof}

We next state the spectral identities for the integral operator
$T_U$; see Lov\'asz~\cite[Chapter~7]{Lovasz2012}.

\begin{lemma}\label{evt:lem:spectral-identities}
Let $(T_Uf)(x)=\int U(x,y)f(y)\dd y$ on real $L^2[0,1]$.
Choose real orthonormal eigenfunctions $f_i$ corresponding to all
nonzero eigenvalues $\lambda_i$ of $T_U$, counted with multiplicity.
Then
\begin{equation}\label{evt:eq:eigenvalue-degree}
 \abs{\lambda_i}\le\int d_U(x)f_i(x)^2\dd x\le\Delta_U,
 \qquad
 \sum_i\lambda_i^2=\int U(x,y)^2\dd x\dd y.
\end{equation}
For every integer $s\ge3$ and every integer $j\ge1$,
\begin{align}
 t(C_s,U)&=\sum_i\lambda_i^s,
 \label{evt:eq:cycle-trace}\\
 \sum_i\lambda_i^2\int d_U(x)^j f_i(x)^2\dd x
 &=\int d_U(x)^jU(x,y)^2\dd x\dd y.
 \label{evt:eq:weighted-HS}
\end{align}
All these series converge absolutely.
\end{lemma}

\begin{proof}
The operator $T_U$ is compact, self-adjoint, and Hilbert--Schmidt.
The spectral theorem gives the eigenfunctions $f_i$; adjoining an
orthonormal basis of $\ker T_U$ gives a complete eigenbasis.
For every real $f\in L^2[0,1]$, symmetry gives
\begin{equation}\label{evt:eq:laplacian-forms}
 \int d_U(x)f(x)^2\dd x\pm\inner{f}{T_Uf}
 =\frac12\int U(x,y)(f(x)\pm f(y))^2\dd x\dd y\ge0.
\end{equation}
Taking $f=f_i$ proves the eigenvalue bound. The formula for
$\sum_i\lambda_i^2$ is the Hilbert--Schmidt norm identity.

For $s\ge3$, the operators $T_U$ and $T_U^{s-1}$ are both
Hilbert--Schmidt, and $U^{\circ(s-1)}$ is the kernel of $T_U^{s-1}$.
Computing their inner product from the kernels and the eigenbasis gives
\[
 t(C_s,U)
 =\int U(x,y)U^{\circ(s-1)}(x,y)\dd x\dd y
 =\inner{T_U}{T_U^{s-1}}_{\mathrm{HS}}
 =\sum_i\lambda_i^s.
\]
For almost every $x$, Parseval's identity also gives
\[
 \int U(x,y)^2\dd y
 =\sum_i\left(\int U(x,y)f_i(y)\dd y\right)^2
 =\sum_i\lambda_i^2 f_i(x)^2.
\]
Multiplying by $d_U(x)^j$ and integrating proves
\eqref{evt:eq:weighted-HS} by Tonelli's theorem. Absolute convergence
follows from $|\lambda_i|\le\Delta_U\le1$ and
$\sum_i\lambda_i^2<\infty$.
\end{proof}

\begin{lemma}\label{evt:lem:nonnegative-decomposition}
With the eigenfunctions and eigenvalues above, for every integer $s\ge3$,
\begin{align}
 M_{s-1}(U)-t(C_s,U)
 ={}&\int d_U(x)^{s-2}U(x,y)(1-U(x,y))\dd x\dd y\notag\\
 &+\sum_i\lambda_i^2
       \left(\int d_U(x)^{s-2}f_i(x)^2\dd x-\lambda_i^{s-2}\right).
 \label{evt:eq:nonnegative-decomposition}
\end{align}
The integral and each summand are nonnegative, and the series
converges absolutely.
\end{lemma}

\begin{proof}
By \eqref{evt:eq:weighted-HS},
\begin{align*}
 M_{s-1}(U)
 &=\int d_U(x)^{s-2}U(x,y)\dd x\dd y\\
 &=\int d_U(x)^{s-2}U(x,y)(1-U(x,y))\dd x\dd y\\
 &\qquad+\sum_i\lambda_i^2\int d_U(x)^{s-2}f_i(x)^2\dd x.
\end{align*}
Subtracting \eqref{evt:eq:cycle-trace} gives the decomposition.
Since $\int f_i(x)^2\dd x=1$, Jensen's inequality and
\eqref{evt:eq:eigenvalue-degree} give, for $j\ge1$,
\[
 \int d_U(x)^j f_i(x)^2\dd x
 \ge\left(\int d_U(x)f_i(x)^2\dd x\right)^j
 \ge|\lambda_i|^j\ge\lambda_i^j.
\]
Thus the integral and each summand are nonnegative. For absolute
convergence, we have
\begin{align*}
 &\sum_i\lambda_i^2
       \left|\int d_U(x)^{s-2}f_i(x)^2\dd x-\lambda_i^{s-2}\right|\\
 &\quad\le\sum_i\lambda_i^2
       \left(\int d_U(x)^{s-2}f_i(x)^2\dd x+|\lambda_i|^{s-2}\right)\\
 &\quad\le2\Delta_U^{s-2}\sum_i\lambda_i^2
       =2\Delta_U^{s-2}\int U(x,y)^2\dd x\dd y<\infty,
\end{align*}
where we used $\int f_i(x)^2\dd x=1$, $d_U\le\Delta_U$, and
\eqref{evt:eq:eigenvalue-degree}.
\end{proof}

\begin{proof}[Proof of \cref{evt:thm:cycle-comparison}]
If $\Delta_U=0$, then $U$ is regular of degree zero, so $U=0$
almost everywhere and both sides of \eqref{evt:eq:cycle-comparison}
are zero. We may therefore assume $\Delta_U>0$.

Fix $i$ with $\lambda_i\ge0$. Since $\int f_i(x)^2\dd x=1$,
let $X$ have probability density $f_i(x)^2$ on $[0,1]$, and set
$Z=d_U(X)$. Then $0\le Z\le\Delta_U$ almost surely, and
\[
 \E Z^j=\int d_U(x)^j f_i(x)^2\dd x\qquad(j\ge1).
\]
Moreover, \eqref{evt:eq:eigenvalue-degree} gives
\[
 0\le\lambda_i\le\int d_U(x)f_i(x)^2\dd x=\E Z.
\]
Since $2\le g-2\le s-2$, we can apply \cref{evt:lem:scalar} with
\[
 a=g-2,\qquad b=s-2,\qquad\alpha=\lambda_i,\qquad\Delta=\Delta_U.
\]
This gives
\begin{equation}\label{evt:eq:positive-eigenvalues}
 \begin{aligned}
 &\int d_U(x)^{s-2}f_i(x)^2\dd x-\lambda_i^{s-2}\\
 &\quad\le\frac{(s-2)(s-3)}{(g-2)(g-3)}\Delta_U^{s-g}
       \left(\int d_U(x)^{g-2}f_i(x)^2\dd x-\lambda_i^{g-2}\right).
 \end{aligned}
\end{equation}
For $\lambda_i<0$, use that $g-2$ is odd to obtain
\begin{align}
 &\int d_U(x)^{s-2}f_i(x)^2\dd x-\lambda_i^{s-2}\notag\\
 &\quad\le\int d_U(x)^{s-2}f_i(x)^2\dd x+|\lambda_i|^{s-2}\notag\\
 &\quad\le\Delta_U^{s-g}
       \left(\int d_U(x)^{g-2}f_i(x)^2\dd x+|\lambda_i|^{g-2}\right)\notag\\
 &\quad=\Delta_U^{s-g}
       \left(\int d_U(x)^{g-2}f_i(x)^2\dd x-\lambda_i^{g-2}\right).
 \label{evt:eq:negative-eigenvalues}
\end{align}
By \eqref{evt:eq:nonnegative-decomposition},
\begin{align*}
 &M_{s-1}(U)-t(C_s,U)\\
 &\quad=\int d_U(x)^{s-2}U(x,y)(1-U(x,y))\dd x\dd y\\
 &\qquad\quad+\sum_i\lambda_i^2
       \left(\int d_U(x)^{s-2}f_i(x)^2\dd x-\lambda_i^{s-2}\right)\\
 &\quad\le\Delta_U^{s-g}
       \int d_U(x)^{g-2}U(x,y)(1-U(x,y))\dd x\dd y\\
 &\qquad\quad+\frac{(s-2)(s-3)}{(g-2)(g-3)}\Delta_U^{s-g}
       \sum_{i:\lambda_i\ge0}\lambda_i^2
       \left(\int d_U(x)^{g-2}f_i(x)^2\dd x-\lambda_i^{g-2}\right)\\
 &\qquad\quad+\Delta_U^{s-g}\sum_{i:\lambda_i<0}\lambda_i^2
       \left(\int d_U(x)^{g-2}f_i(x)^2\dd x-\lambda_i^{g-2}\right)\\
 &\quad\le\frac{(s-2)(s-3)}{(g-2)(g-3)}\Delta_U^{s-g}
       \Biggl(\int d_U(x)^{g-2}U(x,y)(1-U(x,y))\dd x\dd y\\
 &\qquad\qquad\qquad+\sum_i\lambda_i^2
       \left(\int d_U(x)^{g-2}f_i(x)^2\dd x-\lambda_i^{g-2}\right)\Biggr)\\
 &\quad=\frac{(s-2)(s-3)}{(g-2)(g-3)}\Delta_U^{s-g}
       \bigl(M_{g-1}(U)-t(C_g,U)\bigr).
\end{align*}
The first inequality follows from \eqref{evt:eq:positive-eigenvalues},
\eqref{evt:eq:negative-eigenvalues}, and
$d_U(x)^{s-2}\le\Delta_U^{s-g}d_U(x)^{g-2}$.
The second follows from $(s-2)(s-3)/((g-2)(g-3))\ge1$ and the
nonnegativity of the integral and summands. The last equality is
\eqref{evt:eq:nonnegative-decomposition} with $s=g$.
\end{proof}

\subsection{The bound for small maximum degree}\label{evt:sec:polynomial}

For $s\ge3$, let $c_s(H)$ be the number of simple $s$-cycles in $H$,
counted without a choice of orientation or starting vertex.

In this subsection, $H$ has girth $g$, where $g$ is odd. Thus $c_g(H)\ge1$ and
$c_s(H)=0$ for $s<g$.

\begin{lemma}\label{evt:lem:cycle-count}
For $g\le s\le v(H)$,
\begin{equation}\label{evt:eq:cycle-count}
 2s\,c_s(H)\le v(H)^{\lceil2s/(g-1)\rceil}\le v(H)^{s-g+3}.
\end{equation}
In particular, $c_g(H)\le v(H)^3$.
\end{lemma}

\begin{proof}
Set $h=(g-1)/2$ and $k=\lceil s/h\rceil$.
Two distinct paths $P,Q$ of length at most $h$ with the same
endpoints would contain a cycle in their union, giving
\[
 g\le e(P\cup Q)\le e(P)+e(Q)\le2h=g-1.
\]
Since $1\le s-(k-1)h\le h$, the map
\[
 (v_0,\ldots,v_{s-1})\longmapsto
 (v_0,v_h,\ldots,v_{(k-1)h})
\]
is injective on oriented $s$-cycles with a specified starting vertex.
Consequently,
\[
 2s\,c_s(H)\le\lvert V(H)^k\rvert=v(H)^k.
\]
Finally, $s\ge g=2h+1$ and $h\ge1$ give
\[
 k=\left\lceil\frac{2s}{g-1}\right\rceil
  =2+\left\lceil\frac{s-g+1}{h}\right\rceil
  \le2+(s-g+1)=s-g+3.
\]
\end{proof}

Deleting an edge of $C_s$ and applying \cref{evt:lem:tree} gives
\begin{equation}\label{evt:eq:Q-below-S}
 0\le t(P_s,U)-t(C_s,U)\le M_{s-1}(U)-t(C_s,U)
 \qquad(s\ge3).
\end{equation}

For $g\le s\le v(H)$, we obtain
\begin{equation}\label{evt:eq:uniform-cycle}
 t(P_s,U)-t(C_s,U)\le
 \begin{cases}
 M_{g-1}(U)-t(C_g,U),&s=g,\\
 v(H)^2\Delta_U^{s-g}\bigl(M_{g-1}(U)-t(C_g,U)\bigr),&g<s\le v(H).
 \end{cases}
\end{equation}
The case $s=g$ follows from \eqref{evt:eq:Q-below-S}.
For $g=3$ and $3<s\le v(H)$,
\begin{align*}
 t(P_s,U)-t(C_s,U)
 &\le(s-2)\Delta_U^{s-3}\bigl(M_2(U)-t(C_3,U)\bigr)\\
 &\le v(H)^2\Delta_U^{s-3}\bigl(M_2(U)-t(C_3,U)\bigr),
\end{align*}
where the first inequality is \cref{evt:lem:triangle-path} and the
last follows from $s-2\le v(H)^2$.
For odd $g\ge5$ and $g<s\le v(H)$,
\begin{align*}
 t(P_s,U)-t(C_s,U)
 &\le M_{s-1}(U)-t(C_s,U)\\
 &\le\frac{(s-2)(s-3)}{(g-2)(g-3)}\Delta_U^{s-g}
       \bigl(M_{g-1}(U)-t(C_g,U)\bigr)\\
 &\le v(H)^2\Delta_U^{s-g}\bigl(M_{g-1}(U)-t(C_g,U)\bigr),
\end{align*}
where the first two inequalities follow from \eqref{evt:eq:Q-below-S}
and \cref{evt:thm:cycle-comparison}, respectively, and the last follows
from $(s-2)(s-3)/((g-2)(g-3))\le s^2\le v(H)^2$.

\begin{proposition}\label{evt:prop:polynomial-local}
Let $H$ have girth $g$, where $g$ is odd. If $e(H)\Delta_U\le1/2$ and
$v(H)\Delta_U\le1/2$, then
\begin{equation}\label{evt:eq:polynomial-local}
 \begin{aligned}
 &t(H,1-U)-\Phi_H(1-q)\\
 &\quad\ge\left(\sum_{v\in V(H)}\binom{d_H(v)}2-c_3(H)-2v(H)^6\Delta_U\right)\operatorname{Var}(d_U)\\
 &\qquad\quad+c_g(H)(1-5v(H)^6\Delta_U)\bigl(M_{g-1}(U)-t(C_g,U)\bigr).
 \end{aligned}
\end{equation}
In particular, $\Delta_U\le(6v(H)^6)^{-1}$ implies
$t(H,1-U)\ge\Phi_H(1-q)$.
\end{proposition}

\begin{proof}
Each triangle contains three adjacent edge pairs, and an adjacent
pair belongs to at most one triangle. Since $H$ contains a cycle,
\[
 \sum_{v\in V(H)}\binom{d_H(v)}2-c_3(H)\ge1.
\]
Assume $\Delta_U>0$.

For forests $F$ with at most two edges,
\[
 t(F,U)=
 \begin{cases}
  1,&|F|=0,\\
  q,&|F|=1,\\
  q^2,&F\text{ has two disjoint edges},\\
  \displaystyle\int d_U(x)^2\dd x=M_2(U),
       &F\text{ has two adjacent edges}.
 \end{cases}
\]
Their contribution to the first sum of \eqref{evt:eq:forest-gamma} is
\begin{align*}
 &\sum_{\substack{F\text{ forest},\,B(F)=\varnothing\\|F|\le2}}
       (-1)^{|F|}\bigl(t(F,U)-q^{|F|}\bigr)\\
 &\quad=\left(\sum_{v\in V(H)}\binom{d_H(v)}2-c_3(H)\right)
       \bigl(M_2(U)-q^2\bigr)\\
 &\quad=\left(\sum_{v\in V(H)}\binom{d_H(v)}2-c_3(H)\right)
       \operatorname{Var}(d_U).
\end{align*}
Here $\sum_{v\in V(H)}\binom{d_H(v)}2$ counts the pairs of adjacent
edges, and exactly $c_3(H)$ of these have nonempty $B(F)$, one for
each triangle with its largest edge removed.

If $B(F)\ne\varnothing$, choose $e=uv\in B(F)$. The edge $e$ and
the $u$--$v$ path in $F$ form a cycle, so $|F|\ge g-1$.
If $|F|=g-1$, this cycle has length $g$, and $F$ is obtained by
deleting its largest edge $e$. Any other edge in $B(F)$ would join
two vertices at distance at most $g-2$ along the path and create a
cycle of length at most $g-1$. Thus $B(F)=\{e\}$, and
\begin{align*}
 I_F(U)
 &=\int\biggl(\prod_{ab\in F}U(x_a,x_b)\biggr)
       (1-U(x_u,x_v))\prod_{w\in V(H)}\dd x_w\\
 &=t(P_g,U)-t(C_g,U).
\end{align*}
Each $g$-cycle gives exactly one such forest by deleting its largest
edge. Hence
\begin{align*}
 &\sum_{\substack{F\text{ forest},\,B(F)\ne\varnothing\\|F|=g-1}}
       (-1)^{|F|}I_F(U)\\
 &\quad=(-1)^{g-1}c_g(H)\bigl(t(P_g,U)-t(C_g,U)\bigr)\\
 &\quad=c_g(H)\bigl(t(P_g,U)-t(C_g,U)\bigr),
\end{align*}
where the last equality follows from $g$ being odd.
Combining these contributions in \eqref{evt:eq:forest-gamma} and
bounding the remaining terms in absolute value gives
\begin{equation}\label{evt:eq:forest-leading}
 \begin{aligned}
 &t(H,1-U)-\Phi_H(1-q)\\
 &\quad\ge\left(\sum_{v\in V(H)}\binom{d_H(v)}2-c_3(H)\right)\operatorname{Var}(d_U)
       +c_g(H)\bigl(t(P_g,U)-t(C_g,U)\bigr)\\
 &\qquad\quad-\sum_{\substack{F\text{ forest},\,B(F)=\varnothing\\|F|\ge3}}
                 \bigl(t(F,U)-q^{|F|}\bigr)
       -\sum_{\substack{F\text{ forest},\,B(F)\ne\varnothing\\|F|\ge g}}I_F(U).
 \end{aligned}
\end{equation}

By \eqref{evt:eq:forest-variance},
\begin{align}
 &\sum_{\substack{F\text{ forest},\,B(F)=\varnothing\\|F|\ge3}}
                 \bigl(t(F,U)-q^{|F|}\bigr)\notag\\
 &\quad\le \operatorname{Var}(d_U)\sum_{r=3}^{e(H)}\binom{e(H)}{r}\binom r2\Delta_U^{r-2}\notag\\
 &\quad=\binom{e(H)}{2}\bigl((1+\Delta_U)^{e(H)-2}-1\bigr)\operatorname{Var}(d_U)
 \le2e(H)^3\Delta_U \operatorname{Var}(d_U).
 \label{evt:eq:polynomial-EX}
\end{align}
Here $e(H)\Delta_U\le1/2$ and comparison with a geometric series give
\[
 (1+\Delta_U)^{e(H)}\le2,\qquad
 (1+\Delta_U)^{e(H)}-1\le2e(H)\Delta_U.
\]

For each $e=uv\in B(F)$, \cref{evt:lem:spanning} gives
\[
 I_F(U)\le\int\prod_{ab\in F}U(x_a,x_b)
                    (1-U(x_u,x_v))\prod_{w\in V(H)}\dd x_w
 \le\Delta_U^j\bigl(t(P_s,U)-t(C_s,U)\bigr),
\]
where $s$ is the length of the cycle formed by $e$ and its path in
$F$, and $j=|F|-s+1$. Here $e$ is the largest edge of the cycle,
and there are at most $\binom{e(H)}{j}$ choices for the remaining edges
of $F$. For $|F|\ge g$, we have $j\ge1$ when $s=g$ and $j\ge0$
when $s>g$. Summing these bounds over all such forests $F$ and all
edges $e\in B(F)$ gives
\begin{equation}\label{evt:eq:polynomial-EC-raw}
 \begin{aligned}
 &\sum_{\substack{F\text{ forest},\,B(F)\ne\varnothing\\|F|\ge g}}I_F(U)\\
 &\quad\le c_g(H)\bigl((1+\Delta_U)^{e(H)}-1\bigr)
                         \bigl(t(P_g,U)-t(C_g,U)\bigr)\\
 &\qquad\quad+(1+\Delta_U)^{e(H)}
          \sum_{s=g+1}^{v(H)}c_s(H)\bigl(t(P_s,U)-t(C_s,U)\bigr).
 \end{aligned}
\end{equation}
By \eqref{evt:eq:cycle-count} and \eqref{evt:eq:uniform-cycle},
\begin{align*}
 &\sum_{s=g+1}^{v(H)}c_s(H)\bigl(t(P_s,U)-t(C_s,U)\bigr)\\
 &\quad\le v(H)^5\bigl(M_{g-1}(U)-t(C_g,U)\bigr)
                        \sum_{j=1}^{v(H)-g}(v(H)\Delta_U)^j\\
 &\quad\le2v(H)^6\Delta_U\bigl(M_{g-1}(U)-t(C_g,U)\bigr).
\end{align*}
It follows that
\begin{align}
 &\sum_{\substack{F\text{ forest},\,B(F)\ne\varnothing\\|F|\ge g}}I_F(U)\notag\\
 &\quad\le(2c_g(H)e(H)+4v(H)^6)\Delta_U
                             \bigl(M_{g-1}(U)-t(C_g,U)\bigr)\notag\\
 &\quad\le5c_g(H)v(H)^6\Delta_U\bigl(M_{g-1}(U)-t(C_g,U)\bigr),
 \label{evt:eq:polynomial-EC}
\end{align}
where $c_g(H)\ge1$ and $2e(H)\le v(H)^2\le v(H)^6$.

For $g=3$, we have $t(P_g,U)=M_{g-1}(U)$. For $g\ge5$,
\cref{evt:lem:tree} and \eqref{evt:eq:moment-variance} give
\begin{align*}
 0\le M_{g-1}(U)-t(P_g,U)
 &\le M_{g-1}(U)-q^{g-1}\\
 &\le\binom{g-1}{2}\Delta_U^{g-3}\operatorname{Var}(d_U).
\end{align*}
Using $c_g(H)\le v(H)^3$ and $\Delta_U\le1$, both cases give
\begin{align*}
 c_g(H)\bigl(t(P_g,U)-t(C_g,U)\bigr)
 &\ge c_g(H)\bigl(M_{g-1}(U)-t(C_g,U)\bigr)
       -\frac{v(H)^5}{2}\Delta_U^2\operatorname{Var}(d_U).
\end{align*}
Combining this with \eqref{evt:eq:forest-leading},
\eqref{evt:eq:polynomial-EX}, and \eqref{evt:eq:polynomial-EC} gives
\eqref{evt:eq:polynomial-local}, since
\[
 2e(H)^3\Delta_U+\frac{v(H)^5}{2}\Delta_U^2
 \le\left(\frac{v(H)^6}{4}+\frac{v(H)^5}{2}\right)\Delta_U
 \le2v(H)^6\Delta_U.
\]
Finally, $\Delta_U\le(6v(H)^6)^{-1}$ satisfies both hypotheses and gives
\[
 \begin{aligned}
 &t(H,1-U)-\Phi_H(1-q)\\
 &\quad\ge\frac23\operatorname{Var}(d_U)
       +\frac{c_g(H)}6\bigl(M_{g-1}(U)-t(C_g,U)\bigr)\ge0.
 \end{aligned}
\]
\end{proof}

\subsection{From maximum degree to edge density}\label{evt:sec:exceptional}

For a graphon $V$, define
\[
 D_H(V)=t(H,1-V)-1+e(H)\,t(K_2,V).
\]
We write $(a)_+=\max\{a,0\}$.

\begin{lemma}\label{evt:lem:missing-edges}
For every graphon $V$, $D_H(V)\ge0$. If $V_1,V_2$ are graphons with
$V_1+V_2\le1$, then
\begin{equation}\label{evt:eq:D-superadditive}
 D_H(V_1+V_2)\ge D_H(V_1)+D_H(V_2).
\end{equation}
\end{lemma}

\begin{proof}
Choose independent uniform positions $(Z_v)_{v\in V(H)}$ in $[0,1]$.
Conditional on these positions, select each edge $uv$ independently
with probability $V(Z_u,Z_v)$, and let $N$ be the number of selected
edges.
Then $\E N=e(H)\,t(K_2,V)$ and $\Prob(N=0)=t(H,1-V)$, so
\begin{equation}\label{evt:eq:D-probability}
 \begin{aligned}
 D_H(V)&=\E N-1+\Prob(N=0)\\
 &=\E\bigl(N-1+\ind_{\{N=0\}}\bigr)=\E(N-1)_+\ge0.
 \end{aligned}
\end{equation}
To prove \eqref{evt:eq:D-superadditive}, assign each edge $uv$ one
of three outcomes
with probabilities $V_1(Z_u,Z_v)$, $V_2(Z_u,Z_v)$, and
$1-V_1(Z_u,Z_v)-V_2(Z_u,Z_v)$, independently conditional on the
positions. Let $N_1$ and $N_2$ be the numbers of edges with the first
and second outcomes, respectively. Then
\[
 (N_1+N_2-1)_+\ge(N_1-1)_++(N_2-1)_+.
\]
Taking expectations proves \eqref{evt:eq:D-superadditive}.
\end{proof}

\begin{lemma}\label{evt:lem:polynomial-remainder}
For $e(H)\ge1$ and $q_1,q_2\ge0$ with $q=q_1+q_2\le1/(2e(H))$,
\begin{equation}\label{evt:eq:polynomial-remainder}
 \left|\Phi_H(1-q)-\Phi_H(1-q_2)+e(H)q_1\right|
 \le2e(H)^2qq_1.
\end{equation}
\end{lemma}

\begin{proof}
If $e(H)=1$, then $\Phi_H(p)=p$ and the left-hand side is zero.
For $e(H)\ge2$, \eqref{evt:eq:forest-chromatic} gives
$\Phi_H'(1)=e(H)$ and shows that the coefficient of $q^j$ in
$\Phi_H(1-q)$ has absolute value at most $\binom{e(H)}{j}$.
Hence, for $0\le x\le1/(2e(H))$,
\begin{align*}
 |\Phi_H''(1-x)|
 &\le\sum_{j=2}^{e(H)}\binom{e(H)}{j} j(j-1)x^{j-2}\\
 &=e(H)(e(H)-1)(1+x)^{e(H)-2}\le2e(H)^2.
\end{align*}
Integrating gives
\[
 \left|e(H)-\Phi_H'(1-x)\right|
 =\left|\int_0^x\Phi_H''(1-y)\dd y\right|
 \le2e(H)^2x.
\]
Since $q-q_2=q_1$, a second integration gives
\begin{align*}
 &\left|\Phi_H(1-q)-\Phi_H(1-q_2)+e(H)q_1\right|\\
 &\quad=\left|\int_{q_2}^q\bigl(e(H)-\Phi_H'(1-x)\bigr)\dd x\right|\\
 &\quad\le\int_{q_2}^q2e(H)^2x\dd x\\
 &\quad=e(H)^2(q^2-q_2^2)=e(H)^2q_1(q+q_2)\le2e(H)^2qq_1.
\end{align*}
\end{proof}

\begin{proposition}\label{evt:prop:polynomial-reduction}
Suppose $H$ has a vertex of degree at least two, and let
$0<\delta\le1/v(H)$. Assume that every graphon $V$ with
$\Delta_V\le\delta$ satisfies
\[
 t(H,1-V)\ge\Phi_H(1-t(K_2,V)).
\]
Then $t(H,1-U)\ge\Phi_H(1-q)$ whenever
\begin{equation}\label{evt:eq:polynomial-reduction}
 q=t(K_2,U)<
 \min\left\{\frac{\delta^2}{2},\frac{\delta}{32e(H)^2}\right\}.
\end{equation}
If $16e(H)^2\delta\le1$, it suffices that $q<\delta^2/2$.
\end{proposition}

\begin{proof}
Let $A=\{x\in[0,1]:d_U(x)>\delta\}$ and let $s=|A|$ be its
size. Decompose $U$ as
\[
 U_1(x,y)=U(x,y)\ind_{\{x\in A\text{ or }y\in A\}},
 \qquad U_2=U-U_1,
\]
and write $q_i=t(K_2,U_i)$. Then $q_1+q_2=q$ and
\[
 s\le\frac q\delta\le\frac\delta2\le\frac1{2v(H)}.
\]
Also $\Delta_{U_2}\le\delta$, so
$t(H,1-U_2)\ge\Phi_H(1-q_2)$ by hypothesis.
For $x\in A$, the average of $U(x,y)$ over $A^c$ satisfies
\begin{equation}\label{evt:eq:h-lower}
 \frac1{1-s}\int_{A^c}U(x,y)\dd y
 \ge\frac{d_U(x)-s}{1-s}
 \ge d_U(x)-s\ge\frac{d_U(x)}2.
\end{equation}
Moreover,
\begin{equation}\label{evt:eq:exceptional-mass}
 q_1=2\int_A d_U(x)\dd x-\int_{A\times A}U(x,y)\dd x\dd y
 \le2\int_A d_U(x)\dd x.
\end{equation}

For a binomial random variable $N$ with parameters $k\in\mathbb Z_{\ge0}$ and
$h\in[0,1]$, we have
\begin{equation}\label{evt:eq:binomial-bound}
 \E(N-1)_+=(1-h)^k-1+kh\ge(k-1)_+h^2.
\end{equation}
The expression is zero for $k=0,1$; for $k\ge2$,
\[
 (1-h)^k-1+kh
 =h\sum_{j=1}^{k-1}\bigl(1-(1-h)^j\bigr)
 \ge h\sum_{j=1}^{k-1}h=(k-1)h^2.
\]
Let $N$ be the number of edges selected according to $U_1$ as in
\eqref{evt:eq:D-probability}. For each $v\in V(H)$, define the event
\[
 E_v=\{Z_v\in A,\ Z_w\in A^c\text{ for every }w\ne v\}.
\]
The events $E_v$ are disjoint, and $(N-1)_+\ge0$, so
\[
 D_H(U_1)=\E(N-1)_+
 \ge\sum_{v\in V(H)}\E\bigl[(N-1)_+\ind_{E_v}\bigr].
\]
On $E_v$, only edges incident to $v$ can be selected. Conditional
on $Z_v=x\in A$ and all other positions lying in $A^c$, $N$ is
binomial with parameters $d_H(v)$ and
$(1-s)^{-1}\int_{A^c}U(x,y)\dd y$.
Since the other positions lie in $A^c$ with probability
$(1-s)^{v(H)-1}$, \eqref{evt:eq:binomial-bound} and
\eqref{evt:eq:h-lower} give
\begin{align*}
 D_H(U_1)
 &\ge(1-s)^{v(H)-1}
       \sum_{v\in V(H)}(d_H(v)-1)_+
       \int_A\left(\frac1{1-s}\int_{A^c}U(x,y)\dd y\right)^2\dd x\\
 &\ge\frac{(1-s)^{v(H)-1}}4
       \left(\sum_{v\in V(H)}(d_H(v)-1)_+\right)
       \int_A d_U(x)^2\dd x\\
 &\ge\frac18\int_A d_U(x)^2\dd x
 \ge\frac\delta8\int_A d_U(x)\dd x
 \ge\frac\delta{16}q_1.
\end{align*}
Here the degree assumption gives $\sum_v(d_H(v)-1)_+\ge1$, and
$(1-s)^{v(H)-1}\ge1-(v(H)-1)s\ge1/2$.

Since $q\le\delta/(32e(H)^2)\le1/(2e(H))$, superadditivity and
\cref{evt:lem:polynomial-remainder} give
\begin{align*}
 t(H,1-U)
 &\ge t(H,1-U_2)+D_H(U_1)-e(H)q_1\\
 &\ge\Phi_H(1-q_2)+D_H(U_1)-e(H)q_1\\
 &\ge\Phi_H(1-q)+D_H(U_1)-2e(H)^2qq_1\\
 &\ge\Phi_H(1-q)+q_1\left(\frac\delta{16}-2e(H)^2q\right)\\
 &\ge\Phi_H(1-q).
\end{align*}
If $16e(H)^2\delta\le1$, then
$\delta^2/2\le\delta/(32e(H)^2)$, so $q<\delta^2/2$ suffices.
\end{proof}

\subsection{The obstruction for even girth}\label{evt:sec:even-girth}

\begin{lemma}\label{evt:lem:minimal-rank}
Suppose $H$ has girth $g$. If $A\subseteq E(H)$ and $H_A$ contains
a cycle, then $r(A)\ge g-1$. Equality holds if and only if deleting
the isolated vertices from $H_A$ leaves a single $g$-cycle.
\end{lemma}

\begin{proof}
A component of $H_A$ containing a cycle has at least $g$ vertices.
A spanning forest of $H_A$ therefore has at least $g-1$ edges, so
$r(A)\ge g-1$. If $r(A)=g-1$,
this component has exactly $g$ vertices and all other vertices are
isolated. Its cycle uses all $g$ vertices, and an additional edge
would be a chord creating a shorter cycle. The converse is immediate.
\end{proof}

\begin{proposition}\label{evt:prop:even-obstruction}
If $H$ has girth $g$, then
\begin{equation}\label{evt:eq:girth-leading-term}
 \Phi_H(1-q)-(1-q)^{e(H)}
 =(-1)^gc_g(H)q^{g-1}+O_H(q^g)
 \quad\text{as }q\downarrow0.
\end{equation}
If $g$ is even, the constant graphon $W\equiv1-1/k$ violates
\eqref{eq:turan-optimal} for every sufficiently large integer $k$.
\end{proposition}

\begin{proof}
Subtracting the binomial expansion of $(1-q)^{e(H)}$ from
\eqref{evt:eq:chromatic-expansion} gives
\[
 \Phi_H(1-q)-(1-q)^{e(H)}
 =\sum_{A\subseteq E(H)}(-1)^{|A|}
             \bigl(q^{r(A)}-q^{|A|}\bigr).
\]
If $H_A$ is a forest, the summand is zero. Each of the $c_g(H)$
shortest cycles
contributes $(-1)^g(q^{g-1}-q^g)$. By \cref{evt:lem:minimal-rank},
every remaining edge set $A$ containing a cycle satisfies $r(A)\ge g$,
proving \eqref{evt:eq:girth-leading-term}.

For even $g$, the leading coefficient is positive, so
$\Phi_H(1-q)>(1-q)^{e(H)}$ for all sufficiently small $q>0$.
For the constant graphon $W\equiv1-q$,
\[
 t(H,W)=(1-q)^{e(H)}<\Phi_H(1-q).
\]
Taking $q=1/k$ and using \eqref{eq:turan-values} shows that the
constant graphon $W\equiv1-1/k$ violates \eqref{eq:turan-optimal}
for every sufficiently large integer $k$.
\end{proof}

\subsection{Proof of the characterisation}\label{evt:sec:classification}

\begin{proof}[Proof of \cref{main:eventual}]
(iii)$\Rightarrow$(ii).
Suppose the girth of $H$ is odd, and set $\delta=(6v(H)^6)^{-1}$.
By \cref{evt:prop:polynomial-local}, every graphon $V$ with
$\Delta_V\le\delta$ satisfies $t(H,1-V)\ge\Phi_H(1-t(K_2,V))$.
The graph $H$ has a vertex of degree at least two, $\delta\le1/v(H)$, and
\[
 16e(H)^2\delta\le16\left(\frac{v(H)^2}{2}\right)^2\frac1{6v(H)^6}
 =\frac2{3v(H)^2}\le1.
\]
Thus \cref{evt:prop:polynomial-reduction} gives
$t(H,1-U)\ge\Phi_H(1-q)$ whenever
\[
 q<\frac{\delta^2}{2}=\frac1{72v(H)^{12}}.
\]
Taking $U=1-W$ proves (ii) with the threshold
\eqref{eq:eventual-polynomial-threshold}.

\smallskip\noindent
(ii)$\Rightarrow$(i).
For every sufficiently large integer $k$ and every graphon $W$ of
edge density $1-1/k$,
\[
 t(H,W)\ge\Phi_H(1-1/k)=t(H,T_k),
\]
where the equality is \eqref{eq:turan-values}. With the threshold in
\eqref{eq:eventual-polynomial-threshold}, it suffices that
$k>72v(H)^{12}$.

\smallskip\noindent
(i)$\Rightarrow$(iii).
Since $H$ is nonbipartite, it contains a cycle. If its girth were
even, \cref{evt:prop:even-obstruction} would contradict (i).
Hence its girth is odd.
\end{proof}

\section{Graphs on at most six vertices}\label{sec:catalogue}

The classification in \cref{main:classification} is obtained by
enumerating the $117$ graphs and proving the bound or constructing
a counterexample for each. The accompanying
\hyperref[par:formalization]{Lean development} contains the complete catalogue, with a
formal proof for every entry. The accompanying notes give the
graph-specific inequalities and certificate constructions.

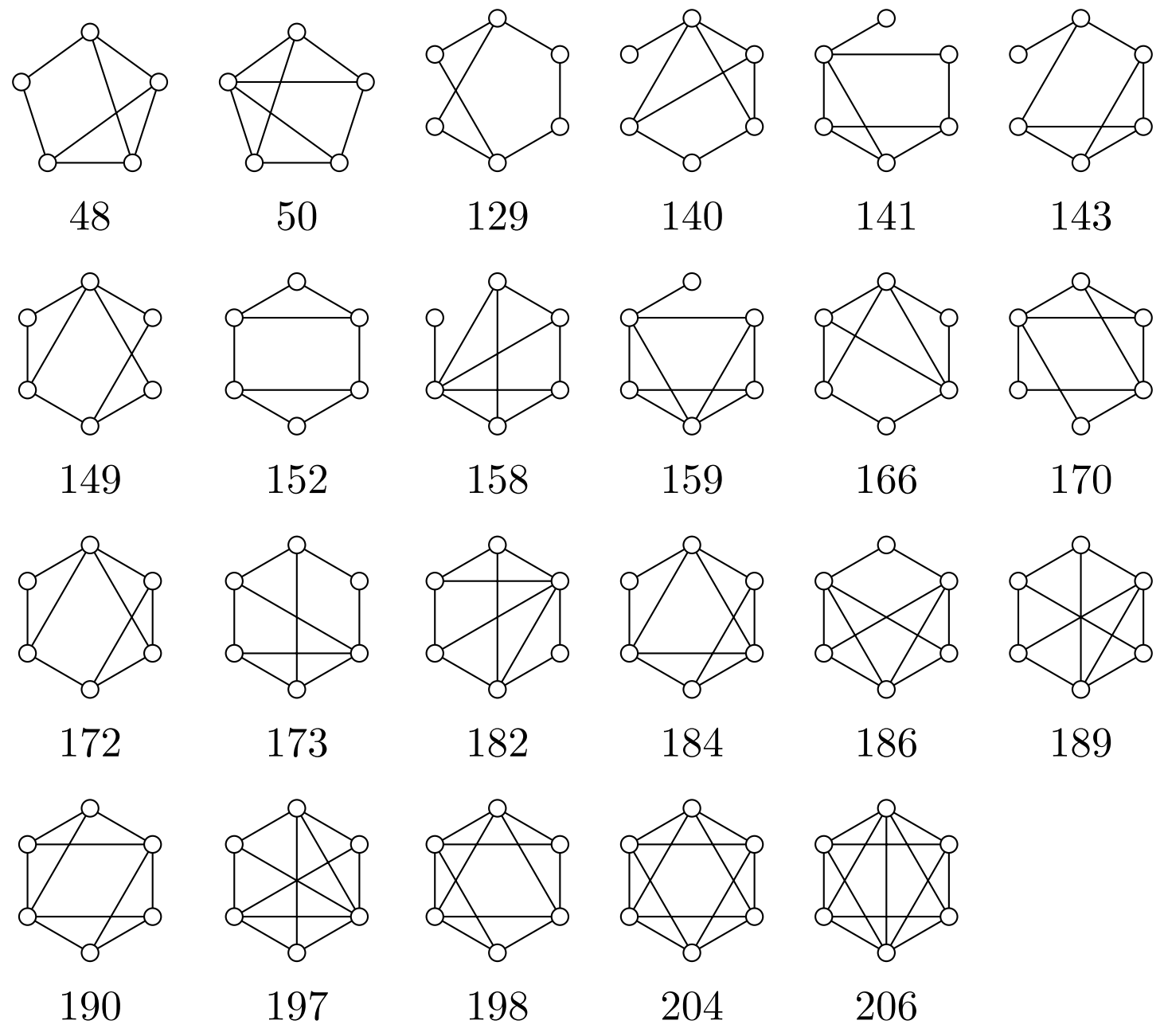
\begin{figure}[htbp]
\centering
\tikzset{gv/.style={circle,draw,fill=white,inner sep=0pt,minimum size=4pt}}
\setlength{\tabcolsep}{6pt}
\begin{tabular}{cccccc}
\begin{tikzpicture}[scale=0.6]\draw (0.000,1.000)--(0.951,0.309) (0.000,1.000)--(-0.951,0.309) (0.000,1.000)--(0.588,-0.809) (-0.588,-0.809)--(0.951,0.309) (-0.588,-0.809)--(-0.951,0.309) (-0.588,-0.809)--(0.588,-0.809) (0.951,0.309)--(0.588,-0.809); \node[gv] at (0.000,1.000) {}; \node[gv] at (0.951,0.309) {}; \node[gv] at (0.588,-0.809) {}; \node[gv] at (-0.588,-0.809) {}; \node[gv] at (-0.951,0.309) {};\end{tikzpicture} &
\begin{tikzpicture}[scale=0.6]\draw (0.000,1.000)--(0.951,0.309) (0.000,1.000)--(-0.588,-0.809) (0.000,1.000)--(-0.951,0.309) (0.951,0.309)--(0.588,-0.809) (0.951,0.309)--(-0.951,0.309) (0.588,-0.809)--(-0.588,-0.809) (0.588,-0.809)--(-0.951,0.309) (-0.588,-0.809)--(-0.951,0.309); \node[gv] at (0.000,1.000) {}; \node[gv] at (0.951,0.309) {}; \node[gv] at (0.588,-0.809) {}; \node[gv] at (-0.588,-0.809) {}; \node[gv] at (-0.951,0.309) {};\end{tikzpicture} &
\begin{tikzpicture}[scale=0.6]\draw (0.000,1.000)--(0.866,0.500) (0.000,1.000)--(-0.866,-0.500) (0.000,1.000)--(-0.866,0.500) (0.866,0.500)--(0.866,-0.500) (0.866,-0.500)--(0.000,-1.000) (0.000,-1.000)--(-0.866,-0.500) (0.000,-1.000)--(-0.866,0.500); \node[gv] at (0.000,1.000) {}; \node[gv] at (0.866,0.500) {}; \node[gv] at (0.866,-0.500) {}; \node[gv] at (0.000,-1.000) {}; \node[gv] at (-0.866,-0.500) {}; \node[gv] at (-0.866,0.500) {};\end{tikzpicture} &
\begin{tikzpicture}[scale=0.6]\draw (0.000,1.000)--(0.866,-0.500) (0.000,1.000)--(-0.866,-0.500) (0.000,1.000)--(0.866,0.500) (0.000,1.000)--(-0.866,0.500) (0.866,-0.500)--(0.866,0.500) (0.866,-0.500)--(0.000,-1.000) (-0.866,-0.500)--(0.866,0.500) (-0.866,-0.500)--(0.000,-1.000); \node[gv] at (0.000,1.000) {}; \node[gv] at (0.866,0.500) {}; \node[gv] at (0.866,-0.500) {}; \node[gv] at (0.000,-1.000) {}; \node[gv] at (-0.866,-0.500) {}; \node[gv] at (-0.866,0.500) {};\end{tikzpicture} &
\begin{tikzpicture}[scale=0.6]\draw (0.000,1.000)--(-0.866,0.500) (0.866,0.500)--(-0.866,0.500) (0.866,0.500)--(0.866,-0.500) (-0.866,0.500)--(0.000,-1.000) (-0.866,0.500)--(-0.866,-0.500) (0.000,-1.000)--(0.866,-0.500) (0.000,-1.000)--(-0.866,-0.500) (0.866,-0.500)--(-0.866,-0.500); \node[gv] at (0.000,1.000) {}; \node[gv] at (0.866,0.500) {}; \node[gv] at (0.866,-0.500) {}; \node[gv] at (0.000,-1.000) {}; \node[gv] at (-0.866,-0.500) {}; \node[gv] at (-0.866,0.500) {};\end{tikzpicture} &
\begin{tikzpicture}[scale=0.6]\draw (0.000,1.000)--(0.866,0.500) (0.000,1.000)--(-0.866,-0.500) (0.000,1.000)--(-0.866,0.500) (0.866,0.500)--(0.866,-0.500) (0.866,0.500)--(0.000,-1.000) (-0.866,-0.500)--(0.866,-0.500) (-0.866,-0.500)--(0.000,-1.000) (0.866,-0.500)--(0.000,-1.000); \node[gv] at (0.000,1.000) {}; \node[gv] at (0.866,0.500) {}; \node[gv] at (0.866,-0.500) {}; \node[gv] at (0.000,-1.000) {}; \node[gv] at (-0.866,-0.500) {}; \node[gv] at (-0.866,0.500) {};\end{tikzpicture} \\
\small 48 & \small 50 & \small 129 & \small 140 & \small 141 & \small 143 \\[6pt]
\begin{tikzpicture}[scale=0.6]\draw (0.000,1.000)--(0.866,0.500) (0.000,1.000)--(-0.866,-0.500) (0.000,1.000)--(-0.866,0.500) (0.000,1.000)--(0.866,-0.500) (0.866,0.500)--(0.000,-1.000) (0.000,-1.000)--(-0.866,-0.500) (0.000,-1.000)--(0.866,-0.500) (-0.866,-0.500)--(-0.866,0.500); \node[gv] at (0.000,1.000) {}; \node[gv] at (0.866,0.500) {}; \node[gv] at (0.866,-0.500) {}; \node[gv] at (0.000,-1.000) {}; \node[gv] at (-0.866,-0.500) {}; \node[gv] at (-0.866,0.500) {};\end{tikzpicture} &
\begin{tikzpicture}[scale=0.6]\draw (0.000,1.000)--(0.866,0.500) (0.000,1.000)--(-0.866,0.500) (0.866,0.500)--(0.866,-0.500) (0.866,0.500)--(-0.866,0.500) (0.866,-0.500)--(0.000,-1.000) (0.866,-0.500)--(-0.866,-0.500) (0.000,-1.000)--(-0.866,-0.500) (-0.866,-0.500)--(-0.866,0.500); \node[gv] at (0.000,1.000) {}; \node[gv] at (0.866,0.500) {}; \node[gv] at (0.866,-0.500) {}; \node[gv] at (0.000,-1.000) {}; \node[gv] at (-0.866,-0.500) {}; \node[gv] at (-0.866,0.500) {};\end{tikzpicture} &
\begin{tikzpicture}[scale=0.6]\draw (0.000,1.000)--(0.866,0.500) (0.000,1.000)--(0.000,-1.000) (0.000,1.000)--(-0.866,-0.500) (0.866,0.500)--(0.866,-0.500) (0.866,0.500)--(-0.866,-0.500) (0.866,-0.500)--(0.000,-1.000) (0.866,-0.500)--(-0.866,-0.500) (0.000,-1.000)--(-0.866,-0.500) (-0.866,-0.500)--(-0.866,0.500); \node[gv] at (0.000,1.000) {}; \node[gv] at (0.866,0.500) {}; \node[gv] at (0.866,-0.500) {}; \node[gv] at (0.000,-1.000) {}; \node[gv] at (-0.866,-0.500) {}; \node[gv] at (-0.866,0.500) {};\end{tikzpicture} &
\begin{tikzpicture}[scale=0.6]\draw (0.000,1.000)--(-0.866,0.500) (0.866,0.500)--(-0.866,0.500) (0.866,0.500)--(0.000,-1.000) (0.866,0.500)--(0.866,-0.500) (-0.866,0.500)--(0.000,-1.000) (-0.866,0.500)--(-0.866,-0.500) (0.000,-1.000)--(0.866,-0.500) (0.000,-1.000)--(-0.866,-0.500) (0.866,-0.500)--(-0.866,-0.500); \node[gv] at (0.000,1.000) {}; \node[gv] at (0.866,0.500) {}; \node[gv] at (0.866,-0.500) {}; \node[gv] at (0.000,-1.000) {}; \node[gv] at (-0.866,-0.500) {}; \node[gv] at (-0.866,0.500) {};\end{tikzpicture} &
\begin{tikzpicture}[scale=0.6]\draw (0.000,1.000)--(-0.866,-0.500) (0.000,1.000)--(0.866,-0.500) (0.000,1.000)--(0.866,0.500) (0.000,1.000)--(-0.866,0.500) (-0.866,-0.500)--(0.000,-1.000) (-0.866,-0.500)--(-0.866,0.500) (0.000,-1.000)--(0.866,-0.500) (0.866,-0.500)--(0.866,0.500) (0.866,-0.500)--(-0.866,0.500); \node[gv] at (0.000,1.000) {}; \node[gv] at (0.866,0.500) {}; \node[gv] at (0.866,-0.500) {}; \node[gv] at (0.000,-1.000) {}; \node[gv] at (-0.866,-0.500) {}; \node[gv] at (-0.866,0.500) {};\end{tikzpicture} &
\begin{tikzpicture}[scale=0.6]\draw (0.000,1.000)--(0.866,0.500) (0.000,1.000)--(0.866,-0.500) (0.000,1.000)--(-0.866,0.500) (0.000,-1.000)--(0.866,-0.500) (0.000,-1.000)--(-0.866,0.500) (0.866,0.500)--(0.866,-0.500) (0.866,0.500)--(-0.866,0.500) (-0.866,-0.500)--(0.866,-0.500) (-0.866,-0.500)--(-0.866,0.500); \node[gv] at (0.000,1.000) {}; \node[gv] at (0.866,0.500) {}; \node[gv] at (0.866,-0.500) {}; \node[gv] at (0.000,-1.000) {}; \node[gv] at (-0.866,-0.500) {}; \node[gv] at (-0.866,0.500) {};\end{tikzpicture} \\
\small 149 & \small 152 & \small 158 & \small 159 & \small 166 & \small 170 \\[6pt]
\begin{tikzpicture}[scale=0.6]\draw (0.000,1.000)--(0.866,0.500) (0.000,1.000)--(-0.866,-0.500) (0.000,1.000)--(-0.866,0.500) (0.000,1.000)--(0.866,-0.500) (0.866,0.500)--(0.000,-1.000) (0.866,0.500)--(0.866,-0.500) (0.000,-1.000)--(-0.866,-0.500) (0.000,-1.000)--(0.866,-0.500) (-0.866,-0.500)--(-0.866,0.500); \node[gv] at (0.000,1.000) {}; \node[gv] at (0.866,0.500) {}; \node[gv] at (0.866,-0.500) {}; \node[gv] at (0.000,-1.000) {}; \node[gv] at (-0.866,-0.500) {}; \node[gv] at (-0.866,0.500) {};\end{tikzpicture} &
\begin{tikzpicture}[scale=0.6]\draw (0.000,1.000)--(0.866,0.500) (0.000,1.000)--(0.000,-1.000) (0.000,1.000)--(-0.866,0.500) (0.866,0.500)--(0.866,-0.500) (0.866,-0.500)--(0.000,-1.000) (0.866,-0.500)--(-0.866,-0.500) (0.866,-0.500)--(-0.866,0.500) (0.000,-1.000)--(-0.866,-0.500) (-0.866,-0.500)--(-0.866,0.500); \node[gv] at (0.000,1.000) {}; \node[gv] at (0.866,0.500) {}; \node[gv] at (0.866,-0.500) {}; \node[gv] at (0.000,-1.000) {}; \node[gv] at (-0.866,-0.500) {}; \node[gv] at (-0.866,0.500) {};\end{tikzpicture} &
\begin{tikzpicture}[scale=0.6]\draw (0.000,1.000)--(0.866,0.500) (0.000,1.000)--(0.000,-1.000) (0.000,1.000)--(-0.866,0.500) (0.866,-0.500)--(0.866,0.500) (0.866,-0.500)--(0.000,-1.000) (0.866,0.500)--(0.000,-1.000) (0.866,0.500)--(-0.866,-0.500) (0.866,0.500)--(-0.866,0.500) (0.000,-1.000)--(-0.866,-0.500) (-0.866,-0.500)--(-0.866,0.500); \node[gv] at (0.000,1.000) {}; \node[gv] at (0.866,0.500) {}; \node[gv] at (0.866,-0.500) {}; \node[gv] at (0.000,-1.000) {}; \node[gv] at (-0.866,-0.500) {}; \node[gv] at (-0.866,0.500) {};\end{tikzpicture} &
\begin{tikzpicture}[scale=0.6]\draw (0.000,1.000)--(0.866,0.500) (0.000,1.000)--(-0.866,-0.500) (0.000,1.000)--(-0.866,0.500) (0.000,1.000)--(0.866,-0.500) (0.866,0.500)--(0.000,-1.000) (0.866,0.500)--(0.866,-0.500) (0.000,-1.000)--(-0.866,-0.500) (0.000,-1.000)--(0.866,-0.500) (-0.866,-0.500)--(-0.866,0.500) (-0.866,-0.500)--(0.866,-0.500); \node[gv] at (0.000,1.000) {}; \node[gv] at (0.866,0.500) {}; \node[gv] at (0.866,-0.500) {}; \node[gv] at (0.000,-1.000) {}; \node[gv] at (-0.866,-0.500) {}; \node[gv] at (-0.866,0.500) {};\end{tikzpicture} &
\begin{tikzpicture}[scale=0.6]\draw (0.000,1.000)--(0.866,0.500) (0.000,1.000)--(-0.866,0.500) (0.866,0.500)--(0.866,-0.500) (0.866,0.500)--(-0.866,-0.500) (0.866,0.500)--(0.000,-1.000) (-0.866,0.500)--(0.866,-0.500) (-0.866,0.500)--(-0.866,-0.500) (-0.866,0.500)--(0.000,-1.000) (0.866,-0.500)--(0.000,-1.000) (-0.866,-0.500)--(0.000,-1.000); \node[gv] at (0.000,1.000) {}; \node[gv] at (0.866,0.500) {}; \node[gv] at (0.866,-0.500) {}; \node[gv] at (0.000,-1.000) {}; \node[gv] at (-0.866,-0.500) {}; \node[gv] at (-0.866,0.500) {};\end{tikzpicture} &
\begin{tikzpicture}[scale=0.6]\draw (0.000,1.000)--(0.866,0.500) (0.000,1.000)--(0.000,-1.000) (0.000,1.000)--(-0.866,0.500) (0.866,-0.500)--(0.866,0.500) (0.866,-0.500)--(0.000,-1.000) (0.866,-0.500)--(-0.866,0.500) (0.866,0.500)--(0.000,-1.000) (0.866,0.500)--(-0.866,-0.500) (0.000,-1.000)--(-0.866,-0.500) (-0.866,-0.500)--(-0.866,0.500); \node[gv] at (0.000,1.000) {}; \node[gv] at (0.866,0.500) {}; \node[gv] at (0.866,-0.500) {}; \node[gv] at (0.000,-1.000) {}; \node[gv] at (-0.866,-0.500) {}; \node[gv] at (-0.866,0.500) {};\end{tikzpicture} \\
\small 172 & \small 173 & \small 182 & \small 184 & \small 186 & \small 189 \\[6pt]
\begin{tikzpicture}[scale=0.6]\draw (0.000,1.000)--(0.866,0.500) (0.000,1.000)--(-0.866,0.500) (0.000,1.000)--(-0.866,-0.500) (0.866,0.500)--(-0.866,0.500) (0.866,0.500)--(0.866,-0.500) (0.866,0.500)--(0.000,-1.000) (-0.866,0.500)--(-0.866,-0.500) (0.866,-0.500)--(0.000,-1.000) (0.866,-0.500)--(-0.866,-0.500) (0.000,-1.000)--(-0.866,-0.500); \node[gv] at (0.000,1.000) {}; \node[gv] at (0.866,0.500) {}; \node[gv] at (0.866,-0.500) {}; \node[gv] at (0.000,-1.000) {}; \node[gv] at (-0.866,-0.500) {}; \node[gv] at (-0.866,0.500) {};\end{tikzpicture} &
\begin{tikzpicture}[scale=0.6]\draw (0.000,1.000)--(0.866,-0.500) (0.000,1.000)--(0.866,0.500) (0.000,1.000)--(0.000,-1.000) (0.000,1.000)--(-0.866,0.500) (0.866,-0.500)--(0.866,0.500) (0.866,-0.500)--(0.000,-1.000) (0.866,-0.500)--(-0.866,-0.500) (0.866,-0.500)--(-0.866,0.500) (0.866,0.500)--(-0.866,-0.500) (0.000,-1.000)--(-0.866,-0.500) (-0.866,-0.500)--(-0.866,0.500); \node[gv] at (0.000,1.000) {}; \node[gv] at (0.866,0.500) {}; \node[gv] at (0.866,-0.500) {}; \node[gv] at (0.000,-1.000) {}; \node[gv] at (-0.866,-0.500) {}; \node[gv] at (-0.866,0.500) {};\end{tikzpicture} &
\begin{tikzpicture}[scale=0.6]\draw (0.000,1.000)--(0.866,0.500) (0.000,1.000)--(-0.866,-0.500) (0.000,1.000)--(0.866,-0.500) (0.000,1.000)--(-0.866,0.500) (0.866,0.500)--(0.866,-0.500) (0.866,0.500)--(-0.866,0.500) (-0.866,-0.500)--(0.000,-1.000) (-0.866,-0.500)--(0.866,-0.500) (-0.866,-0.500)--(-0.866,0.500) (0.000,-1.000)--(0.866,-0.500) (0.000,-1.000)--(-0.866,0.500); \node[gv] at (0.000,1.000) {}; \node[gv] at (0.866,0.500) {}; \node[gv] at (0.866,-0.500) {}; \node[gv] at (0.000,-1.000) {}; \node[gv] at (-0.866,-0.500) {}; \node[gv] at (-0.866,0.500) {};\end{tikzpicture} &
\begin{tikzpicture}[scale=0.6]\draw (0.000,1.000)--(0.866,0.500) (0.000,1.000)--(0.866,-0.500) (0.000,1.000)--(-0.866,-0.500) (0.000,1.000)--(-0.866,0.500) (0.866,0.500)--(0.866,-0.500) (0.866,0.500)--(0.000,-1.000) (0.866,0.500)--(-0.866,0.500) (0.866,-0.500)--(0.000,-1.000) (0.866,-0.500)--(-0.866,-0.500) (0.000,-1.000)--(-0.866,-0.500) (0.000,-1.000)--(-0.866,0.500) (-0.866,-0.500)--(-0.866,0.500); \node[gv] at (0.000,1.000) {}; \node[gv] at (0.866,0.500) {}; \node[gv] at (0.866,-0.500) {}; \node[gv] at (0.000,-1.000) {}; \node[gv] at (-0.866,-0.500) {}; \node[gv] at (-0.866,0.500) {};\end{tikzpicture} &
\begin{tikzpicture}[scale=0.6]\draw (0.000,1.000)--(0.866,0.500) (0.000,1.000)--(0.866,-0.500) (0.000,1.000)--(0.000,-1.000) (0.000,1.000)--(-0.866,-0.500) (0.000,1.000)--(-0.866,0.500) (0.866,0.500)--(0.866,-0.500) (0.866,0.500)--(0.000,-1.000) (0.866,0.500)--(-0.866,0.500) (0.866,-0.500)--(0.000,-1.000) (0.866,-0.500)--(-0.866,-0.500) (0.000,-1.000)--(-0.866,-0.500) (0.000,-1.000)--(-0.866,0.500) (-0.866,-0.500)--(-0.866,0.500); \node[gv] at (0.000,1.000) {}; \node[gv] at (0.866,0.500) {}; \node[gv] at (0.866,-0.500) {}; \node[gv] at (0.000,-1.000) {}; \node[gv] at (-0.866,-0.500) {}; \node[gv] at (-0.866,0.500) {};\end{tikzpicture} &
 \\
\small 190 & \small 197 & \small 198 & \small 204 & \small 206 &  \\[6pt]
\end{tabular}
\caption{The $23$ graphs that violate the Goodman-type lower bound,
labelled by their numbers in the Graph Atlas.}
\label{fig:catalogue-negative}
\end{figure}

For the positive cases, the general results of
\cref{sec:chordal,sec:lifting} are supplemented by estimates for
individual graphs. These estimates express homomorphism densities
as integrals with one or more vertex positions fixed, then apply
the weighted triangle inequality, moment bounds, and
Cauchy--Schwarz. The remaining cases follow from sum-of-squares
identities. For the negative cases, finite step graphons give
strict counterexamples. 
The $23$ violating graphs are shown in \cref{fig:catalogue-negative}.

% BEGIN CLASSIFICATION FIGURE
% END CLASSIFICATION FIGURE

\subsection{Density ranges for the violating graphs}\label{sec:catalogue-refinements}

Of the $23$ graphs in \cref{fig:catalogue-negative}, all except
graph $129$ have odd girth and hence eventually satisfy the
Goodman-type lower bound by \cref{main:eventual}. Graph $129$ has
even girth, so the bound fails for all $p<1$ sufficiently close
to $1$. Thus there is no $p_0<1$ such that the bound holds at every
density $p\in(p_0,1)$ for this graph.
\Cref{tab:catalogue-refinements} gives density ranges on which the
bound holds, together with results on the optimality of $T_k$.

\begin{table}[htbp]
\centering
\small
\renewcommand{\arraystretch}{1.16}
\begin{tabularx}{\linewidth}{@{}>{\raggedright\arraybackslash}p{.30\linewidth}
 >{\raggedright\arraybackslash}X
 >{\raggedright\arraybackslash}p{.30\linewidth}@{}}
\toprule
Graph numbers & Lower bound range & Optimality of $T_k$ \\
\midrule
152, 166, 172, 182, 184
 & $p=2/3$ or $p\ge 3/4$ & Optimal for all $k\ge\chi(H)=3$ \\
\addlinespace
206 & $p\ge 3/4$ & Optimal for all $k\ge\chi(H)=4$ \\
\addlinespace
48, 50, 140, 141, 143, 158, 159, 197, 198
 & $p\ge 3/4$ & Suboptimal at $k=3$; optimal for $k\ge4$ \\
\addlinespace
186 & $p=3/4$ or $p\ge 4/5$ & Suboptimal at $k=3$; optimal for $k\ge4$ \\
\addlinespace
190 & $p\ge 4/5$ & Suboptimal at $k=3,4$; optimal for $k\ge5$ \\
\addlinespace
173 & $p\ge 4/5$ & Suboptimal at $k=3$; optimal for $k\ge5$; $k=4$ open \\
\addlinespace
149, 170 & $p>1-(72\cdot6^{12})^{-1}$ & Suboptimal at $k=3,4,5$;
optimal for $k>72\cdot6^{12}$ \\
\addlinespace
189, 204 & $p>1-(72\cdot6^{12})^{-1}$ & Suboptimal at $k=3,4$;
optimal for $k>72\cdot6^{12}$ \\
\addlinespace
129 & Fails near $p=1$ & Suboptimal at every $k\ge4$ \\
\bottomrule
\end{tabularx}
\caption{Density ranges and Tur\'an optimality for the $23$ graphs
that violate the Goodman-type lower bound.}
\label{tab:catalogue-refinements}
\end{table}

For the $94$ graphs satisfying \eqref{eq:goodman}, $T_k$ is optimal
for every integer $k\ge\chi(H)$. The first two rows of
\cref{tab:catalogue-refinements} give six further graphs with this
property. 

Write $H_{173}$ for graph $173$. At edge density $3/4$, we have
\begin{equation}\label{eq:atlas173-global}
 \frac{14056671}{200000000}
 \le \inf_{\substack{W\text{ a graphon}\\t(K_2,W)=3/4}}t(H_{173},W)
 \le \frac9{128}=t(H_{173},T_4).
\end{equation}
The two bounds differ by $5829/200000000$. Whether $T_4$ attains the
minimum remains open.

\subsection{Constructions of counterexamples}\label{sec:catalogue-counterexamples}

In this subsection, we construct step graphons that violate the
Goodman-type lower bound for each of the $23$ graphs in
\cref{fig:catalogue-negative}.

For graphons $U$ and $V$, their \emph{tensor product} is the graphon
on the product probability space $[0,1]^2$ given by
\[
 (U\otimes V)((x,x'),(y,y'))=U(x,y)V(x',y').
\]
Homomorphism densities are multiplicative under this operation:
\[
 t(H,U\otimes V)=t(H,U)\,t(H,V).
\]

For the graphs $129,149,170,173,186,189,197,198,204$, the graphon
$T_{2k-1}\otimes T_{2k}$, whose edge density is
\[
 \left(1-\frac1{2k-1}\right)\left(1-\frac1{2k}\right)=1-\frac1k,
\]
satisfies $t(H, T_{2k-1} \otimes T_{2k}) < t(H, T_k)$ for each $k$
listed as suboptimal in \cref{tab:catalogue-refinements}.
At $k=3$, the same construction also gives a counterexample for
graph $190$.

For graphs $182$ and $184$, take $W=T_4\otimes T_4$. Then
\[
 t(K_2,W)=\frac9{16},
 \qquad
 t(H,W)=\frac{81}{65536}
 <\frac{351}{262144}=\Phi_H\!\left(\frac9{16}\right).
\]

For graph $129$, the constant graphon $W\equiv1-1/k$ gives
\[
 t(H,T_k)-t(H,W)
 =\frac{(k-1)(k^3-5k^2+6k-1)}{k^7}>0\qquad(k\ge4).
\]
Indeed, substituting $k=u+4$ turns the cubic factor into
$u^3+7u^2+14u+7$.

For graphs $166$ and $172$, partition $[0,1]$ into four equal
parts $A_1,A_2,A_3,A_4$. Define $W$ to be constant on each
$A_i\times A_j$, with value given by the $(i,j)$-entry of the following matrix:
\[
 \begin{pmatrix}
  0 & \frac14 & 1 & 1\\[2pt]
  \frac14 & 0 & 1 & 1\\[2pt]
  1 & 1 & 0 & \frac14\\[2pt]
  1 & 1 & \frac14 & 0
 \end{pmatrix}.
\]
Its edge density is $9/16$, and
\[
 t(H,W)=
 \begin{cases}
  465/262144,&\text{for graph }166,\\
  3/2048,&\text{for graph }172,
 \end{cases}
 \quad <\frac{477}{262144}=\Phi_H\!\left(\frac9{16}\right).
\]

For graph $206$, partition $[0,1]$ into three equal parts
$A_1,A_2,A_3$. Define $W$ to be constant on each $A_i\times A_j$,
with value given by the $(i,j)$-entry of
\[
 \begin{pmatrix}
  \frac14 & 1 & 1\\[2pt]
  1 & 0 & 1\\[2pt]
  1 & 1 & 0
 \end{pmatrix}.
\]
Then
\[
 t(K_2,W)=\frac{25}{36},
 \qquad
 t(H,W)=\frac{77731841}{48922361856}
 <\frac{16975}{10077696}=\Phi_H\!\left(\frac{25}{36}\right).
\]

For graph $190$ at $k=4$, partition $[0,1]$ into four equal parts
$A_1,A_2,A_3,A_4$. For $0\le\varepsilon\le1$, define
$W_\varepsilon$ to be constant on each $A_i\times A_j$, with value
given by the $(i,j)$-entry of the following matrix:
\[
 \begin{pmatrix}
  \varepsilon & 1-\varepsilon/3 & 1-\varepsilon/3 & 1-\varepsilon/3\\[2pt]
  1-\varepsilon/3 & \varepsilon & 1-\varepsilon/3 & 1-\varepsilon/3\\[2pt]
  1-\varepsilon/3 & 1-\varepsilon/3 & \varepsilon & 1-\varepsilon/3\\[2pt]
  1-\varepsilon/3 & 1-\varepsilon/3 & 1-\varepsilon/3 & \varepsilon
 \end{pmatrix}.
\]
The edge density remains $3/4$, while
$\left.\frac{\dd}{\dd\varepsilon}t(H,W_\varepsilon)
\right|_{\varepsilon=0}=-1/256$, so sufficiently small positive
$\varepsilon$ decreases the $H$-density.

For graphs $48,50,140,141,143,158$ at $k=3$, use the following
four-part graphon. Partition
$[0,1]$ into a set $A_1$ of size $\delta=1/2000$ and sets
$A_2,A_3,A_4$ of size $(1-\delta)/3$ each. Define $W$ to be
constant on each product of two parts, with values given by the
following matrix:
\[
 \begin{pmatrix}
  0 & 1 & \frac{23}{50} & \frac{23}{50}\\[2pt]
  1 & \varepsilon & 1 & 1\\[2pt]
  \frac{23}{50} & 1 & \varepsilon & 1\\[2pt]
  \frac{23}{50} & 1 & 1 & \varepsilon
 \end{pmatrix}.
\]
The choice $\varepsilon=8046/99900025$ makes the edge density
exactly $2/3$. For each of the six graphs, $t(H,W)<t(H,T_3)$; the
gaps $t(H,T_3)-t(H,W)$ for graphs $48,50,140,141,143,158$ are
approximately $4.10\times10^{-7}$, $2.24\times10^{-8}$,
$3.30\times10^{-7}$, $1.16\times10^{-7}$, $8.00\times10^{-8}$, and
$1.08\times10^{-7}$, respectively.

For graph $159$, a three-part graphon gives a counterexample at
$k=3$. Partition $[0,1]$ into sets $A_1,A_2,A_3$ of sizes $13/40$,
$7/20$, $13/40$, respectively. Define $W$ to be constant on each
$A_i\times A_j$, with value given by the $(i,j)$-entry of the
following matrix:
\[
 \begin{pmatrix}
  0 & 1 & \frac{49}{50}\\[2pt]
  1 & \frac{557}{14700} & 1\\[2pt]
  \frac{49}{50} & 1 & 0
 \end{pmatrix}.
\]
Its edge density is exactly $2/3$, and
$t(H,W)<t(H,T_3)=4/243$, with $t(H,T_3)-t(H,W)\approx1.744\times10^{-6}$.

For graph $152$, let $0<\varepsilon<1/2$ and partition $[0,1]$
into sets $A_1,A_2,B,S_1,S_2$ of sizes
\[
 \left(\frac12-\varepsilon\right)\frac{624}{625},\quad
 \left(\frac12-\varepsilon\right)\frac1{625},\quad
 \frac12,\quad \frac{\varepsilon}{10},\quad
 \frac{9\varepsilon}{10}.
\]
Define $W_\varepsilon$ to be constant on each product of two parts,
with values given by the following matrix, whose rows and columns
are ordered as $A_1,A_2,B,S_1,S_2$:
\[
 \begin{pmatrix}
  0 & 0 & 1 & 1 & 0\\[2pt]
  0 & 0 & 1 & 0 & 1\\[2pt]
  1 & 1 & 0 & \frac1{16} & 1\\[2pt]
  1 & 0 & \frac1{16} & 0 & 0\\[2pt]
  0 & 1 & 1 & 0 & 0
 \end{pmatrix}.
\]
Let $p_\varepsilon=t(K_2,W_\varepsilon)$. Summing over the partition
classes gives
\begin{align*}
 p_\varepsilon
 &=\frac12+\frac{753}{100000}\varepsilon
          -\frac{633}{3125}\varepsilon^2,\\
 t(H,W_\varepsilon)
 &=\frac{2224881}{80000000000}\varepsilon^2(1-2\varepsilon)^2.
\end{align*}
For this graph, the Goodman-type bound \eqref{eq:goodman}
would require
\[
 t(H,W_\varepsilon)\ge
 p_\varepsilon(2p_\varepsilon-1)^2
 (3p_\varepsilon^2-3p_\varepsilon+1)
 \qquad\text{whenever }p_\varepsilon>\frac12.
\]
The formula for $p_\varepsilon$ shows that $p_\varepsilon>1/2$
for all sufficiently small positive $\varepsilon$. Moreover,
\begin{align*}
 \lim_{\varepsilon\downarrow0}
 \frac{p_\varepsilon(2p_\varepsilon-1)^2
       (3p_\varepsilon^2-3p_\varepsilon+1)}{\varepsilon^2}
 &=\frac12\left(\frac{753}{50000}\right)^2\frac14\\
 &=\frac{2268036}{80000000000}\\
 &>\frac{2224881}{80000000000}
 =\lim_{\varepsilon\downarrow0}
   \frac{t(H,W_\varepsilon)}{\varepsilon^2}.
\end{align*}
Hence, for all sufficiently small positive $\varepsilon$,
\[
 t(H,W_\varepsilon)<
 p_\varepsilon(2p_\varepsilon-1)^2
 (3p_\varepsilon^2-3p_\varepsilon+1),
\]
so $W_\varepsilon$ is a counterexample to \eqref{eq:goodman}.

\section{Concluding remarks}\label{sec:open-problems}

Theorems~\ref{main:pure} and~\ref{main:tp} establish the bound for
two classes of chordal graphs. The following construction shows
that attaching leaves to a clique can produce a chordal graph
that violates the bound.

\begin{proposition}\label{rem:chordal-leaves}
Let $L_{20}$ be obtained from a triangle by attaching twenty
pendant leaves to each of its three vertices. Then $L_{20}$ is a
connected chordal graph with $63$ vertices and $63$ edges that does
not satisfy the Goodman-type lower bound.
\end{proposition}

\begin{proof}
The polynomial $\Phi_{L_{20}}$ is
\[
  \Phi_{L_{20}}(p)=p^{61}(2p-1).
\]
Partition $[0,1]$ into sets $A,B$ of sizes $2/5,3/5$, and define
a graphon $W$ by
\[
 W|_{A\times A}=0,\qquad
 W|_{A\times B}=W|_{B\times A}=1,\qquad
 W|_{B\times B}=\frac19.
\]
Its degree function equals $3/5$ on $A$ and $7/15$ on $B$, so its
edge density is $p=13/25>1/2$. A nonzero contribution to the
triangle density has either one vertex in $A$ or all three
vertices in $B$. Integrating the leaves and summing these two
contributions gives the exact inequality
\[
 t(L_{20},W)
 =\frac6{125}\left(\frac{49}{375}\right)^{20}
  +\frac1{3375}\left(\frac7{15}\right)^{60}
 <\frac1{25}\left(\frac{13}{25}\right)^{61}
 =\Phi_{L_{20}}(13/25).
\]
Thus $W$ violates the Goodman-type lower bound for $L_{20}$.
\end{proof}

\begin{problem}[Chordal classification]\label{prob:chordal}
Characterise the connected nonbipartite chordal graphs that satisfy the
Goodman-type lower bound. Can a necessary and sufficient condition be
expressed in terms of their maximal cliques and the way those cliques
intersect?
\end{problem}

The graphon in the proof of \cref{rem:chordal-leaves} has edge
density $13/25$, which is not a Tur\'an density, so it does not rule
out the optimality of $T_k$ for $L_{20}$. We ask whether $T_k$ is
optimal for every graph in a larger class, which contains every
chordal graph and every odd cycle.

\begin{problem}[Graphs with no induced even cycle]\label{prob:even-hole-free}
Let $H$ be a nonbipartite graph with no induced even cycle. Is $T_k$
optimal for every integer $k\ge\chi(H)$? Equivalently, does
\[
 t(H,W)\ge\frac{P_H(k)}{k^{v(H)}}
\]
hold for every integer $k\ge\chi(H)$ and every graphon $W$ with
$t(K_2,W)=1-1/k$?
\end{problem}

A shortest cycle of such a graph is induced and hence odd, so
\cref{main:eventual} gives the optimality of $T_k$ for all
sufficiently large $k$.

The graphs of \cref{main:classification} also support this question.
Each of the $23$ graphs in \cref{fig:catalogue-negative} contains an
induced $C_4$. Hence all $75$ of the $117$ graphs with no induced even
cycle satisfy the Goodman-type lower bound.

We finally ask for the optimal exponent of $v(H)$ in the density
threshold of \cref{main:eventual}. The following construction shows
that every exponent below $4$ is insufficient.

\begin{proposition}\label{prop:quartic-obstruction}
For every sufficiently large odd integer $n$, there is an
$n$-vertex graph $H$ of girth $3$ and chromatic number $3$ such that
the constant graphon $W\equiv p$, where
\[
 p=1-\frac{128}{(n-1)^4},
\]
satisfies $t(H,W)<\Phi_H(p)$. 
\end{proposition}

\begin{proof}
Let $r=(n-1)/2$ and let $H$ be the graph obtained from
$K_{r,r}$ by adding one vertex adjacent to both ends of a chosen
edge. Denote this graph by $J_r$. The new vertex completes a
triangle and can receive a third colour in a $2$-colouring of
$K_{r,r}$. Thus $J_r$ has girth $3$ and chromatic number $3$.
Moreover, $v(J_r)=2r+1$, $e(J_r)=r^2+2$, and
$P_{J_r}(x)=(x-2)P_{K_{r,r}}(x)$.
Let $\alpha_1,\ldots,\alpha_{2r}$ be the roots of $P_{K_{r,r}}$.
Sokal's bound~\cite{Sokal2001} states that every root of the chromatic
polynomial of a graph has modulus at most eight times its maximum
degree. Since $K_{r,r}$ has maximum degree $r$, we have
$|\alpha_i|\le8r$ for every $i$.
Since $P_{K_{r,r}}$ is monic of degree $2r$,
\[
 \Phi_{K_{r,r}}(1-q)
 =q^{2r}P_{K_{r,r}}(1/q)
 =\prod_{i=1}^{2r}(1-\alpha_iq).
\]
For $0<q\le1/(16r)$, we have $|\alpha_iq|\le1/2$, so the
convergent series $\log(1-z)=-\sum_{j\ge1}z^j/j$ gives
\begin{align*}
 \log\frac{\Phi_{K_{r,r}}(1-q)}{(1-q)^{r^2}}
 &=\sum_{i=1}^{2r}\log(1-\alpha_iq)-r^2\log(1-q)\\
 &=\sum_{j\ge1}\frac{r^2-\sum_{i=1}^{2r}\alpha_i^j}{j}q^j\\
 &=\left(r^2-\sum_{i=1}^{2r}\alpha_i\right)q
   +\frac12\left(r^2-\sum_{i=1}^{2r}\alpha_i^2\right)q^2\\
 &\quad+\frac13\left(r^2-\sum_{i=1}^{2r}\alpha_i^3\right)q^3
   +O(r^5q^4).
\end{align*}
Since $K_{r,r}$ has $r^2$ edges and no triangles,
\eqref{evt:eq:chromatic-expansion} and Vieta's formulas give
\begin{align*}
 \sum_{i=1}^{2r}\alpha_i
 &=\sum_{\substack{A\subseteq E(K_{r,r})\\|A|=1}}1=r^2, \qquad
 \sum_{1\le i<j\le2r}\alpha_i\alpha_j
 =\sum_{\substack{A\subseteq E(K_{r,r})\\|A|=2}}1
 =\binom{r^2}{2}.
\end{align*}
Hence
\begin{align*}
 \sum_{i=1}^{2r}\alpha_i^2
 &=\left(\sum_{i=1}^{2r}\alpha_i\right)^2
   -2\sum_{1\le i<j\le2r}\alpha_i\alpha_j
 =r^4-2\binom{r^2}{2}=r^2.
\end{align*}
Thus the terms with $j=1$ and $j=2$ are zero.

For the cubic term, the edge subsets contributing to the coefficient
of $q^3$ in \eqref{evt:eq:chromatic-expansion} are precisely the
three-edge forests and the four-cycles, with signs $-1$ and $+1$,
respectively. Every three-edge subset of $E(K_{r,r})$ is a forest.
Each four-cycle is obtained by choosing two vertices from each
part, giving $\binom r2^2$ choices. Vieta's formulas therefore give
\[
 \sum_{1\le i<j<k\le2r}\alpha_i\alpha_j\alpha_k
 =\sum_{\substack{A\subseteq E(K_{r,r})\\|A|=3}}1-\binom r2^2
 =\binom{r^2}{3}-\binom r2^2.
\]
It follows that
\begin{align*}
 \sum_{i=1}^{2r}\alpha_i^3
 &=\left(\sum_i\alpha_i\right)^3
   -3\left(\sum_i\alpha_i\right)
     \left(\sum_{i<j}\alpha_i\alpha_j\right)
   +3\sum_{i<j<k}\alpha_i\alpha_j\alpha_k\\
 &=r^6-3r^2\binom{r^2}{2}
   +3\left(\binom{r^2}{3}-\binom r2^2\right)\\
 &=r^2-3\binom r2^2.
\end{align*}
Thus the coefficient of $q^3$ in the logarithmic expansion is
\[
 \frac{r^2-\sum_{i=1}^{2r}\alpha_i^3}{3}=\binom r2^2.
\]
For the remaining terms, $|\sum_i\alpha_i^j|\le2r(8r)^j$ gives
\begin{align*}
 \left|\sum_{j\ge4}\frac{r^2-\sum_{i=1}^{2r}\alpha_i^j}{j}q^j\right|
 &\le\frac{r^2q^4}{4(1-q)}
     +\frac{r(8rq)^4}{2(1-8rq)}\\
 &\le\frac12r^2q^4+r(8rq)^4
 =O(r^5q^4),
\end{align*}
where the last inequality follows from $q\le1/2$ and $8rq\le1/2$.
Consequently,
\begin{align*}
 \log\frac{\Phi_{J_r}(1-q)}{(1-q)^{r^2+2}}
 &=\log\frac{1-2q}{(1-q)^2}
   +\log\frac{\Phi_{K_{r,r}}(1-q)}{(1-q)^{r^2}}\\
 &=-q^2+\left(\binom r2^2-2\right)q^3+O(r^5q^4).
\end{align*}
Taking $q=8/r^4$ gives
\[
 \log\frac{\Phi_{J_r}(p)}{p^{r^2+2}}
 =\frac{64}{r^8}+O(r^{-9})>0
\]
for all sufficiently large $r$. Thus
$t(J_r,W)=p^{r^2+2}<\Phi_{J_r}(p)$.
\end{proof}

Together, \cref{main:eventual,prop:quartic-obstruction} show that
the optimal uniform exponent lies between $4$ and $12$.

\begin{problem}[Optimal exponent in the density threshold]\label{prob:optimal-exponent}
Determine the infimum of the exponents $\gamma$ for which there is
a constant $C>0$ such that
\[
 t(H,W)\ge\Phi_H(p)
 \qquad\text{whenever}\qquad
 p=t(K_2,W)>1-\frac1{Cn^\gamma},
\]
for every integer $n\ge3$, every graph $H$ of odd girth with at most
$n$ vertices, and every graphon $W$.
In particular, does such a bound hold with $\gamma=4$?
\end{problem}

\paragraph{Formal verification.}\label{par:formalization}
The main results in \cref{sec:introduction}, the density ranges and
Tur\'an-optimality statements in \cref{tab:catalogue-refinements},
the counterexamples in \cref{sec:catalogue-counterexamples},
\cref{rem:chordal-leaves}, and the upper bound in
\eqref{eq:atlas173-global} have been formally verified in Lean~4
using Mathlib. The development contains a proof for each of the
$117$ graphs, including the rational sum-of-squares certificates
used in the positive cases and exact counterexamples for the
negative cases. It is available at
\url{https://github.com/mekty2012/GoodmanTypeBound}; the module
\nolinkurl{Goodman/Paper.lean} lists the correspondence between
the paper's statements and their Lean declarations.

The formal proof deviates from the proof in this paper for ease of
formalisation. 
First, in the proof of \cref{evt:lem:forest}, it proves the required absolute-error bound
directly by leaf integration and $L^2$ estimates, instead of using
Sidorenko's tree inequalities. 
Second, in the proof of \cref{evt:thm:cycle-comparison},
it replaces the weighted Parseval identity \eqref{evt:eq:weighted-HS}
by Bessel's inequality and controls an additional nonnegative residual
term in the spectral decomposition. 

\paragraph{AI usage disclosure.}
The initial proofs of
\cref{main:pure,main:odd-join} were generated with GPT-5.5 through
ChatGPT Pro. The author proved the general clique-tree criterion in
\cref{main:profile} and the result for trivially perfect graphs in
\cref{main:tp} by extending these proof ideas. The proof of
\cref{main:eventual} was obtained using GPT-6 Pro, and the $117$
cases in \cref{main:classification} were resolved using GPT-5.6 Codex.
The Lean formalisation was carried out using Claude Opus~5 and~5.5.
Several attempts with ChatGPT~6 Pro produced neither proofs nor
counterexamples for the open problems. 
GPT-6.1 Sol also assisted with drafting proofs and manuscript text.
The author checked and revised this material.

\paragraph{Acknowledgements}
I would like to thank Hongseok Yang for helpful discussions. I was supported by the National Research Foundation of Korea (NRF) grant funded by the Korean Government (MSIT) (No. RS-2023-00279680).

\bibliographystyle{abbrvnat}
\bibliography{paper}
\end{document}